\documentclass[11pt]{amsart}

\usepackage[utf8x]{inputenc}
\usepackage[margin=1.3in]{geometry}
\usepackage{microtype}
\usepackage[T1]{fontenc}
\usepackage{ae,aecompl}
\usepackage{times}
\usepackage{tikz-cd}
\usepackage{mathrsfs} 
\usepackage{harmony}
\usepackage{verbatim,amsmath,amsthm,amsfonts,amssymb,latexsym,graphicx,mathtools,extpfeil,color,mathabx}
\usepackage{tikz}
\usepackage{epstopdf,pinlabel}
\usepackage{graphicx}
\usepackage{booktabs}
\usepackage{caption}
\usepackage[framemethod=TikZ]{mdframed}
\usepackage{tabularx}
\usepackage{subcaption}
\usepackage{slashed}
\DeclareMathAlphabet{\mathpzc}{OT1}{pzc}{m}{it}
\usepackage[colorlinks,pagebackref,hypertexnames=false]{hyperref} \usepackage[alphabetic,backrefs,msc-links]{amsrefs}

\usepackage{aliascnt}
\numberwithin{equation}{section}

\newcommand{\Z}{\mathbb{Z}}

\newcommand{\R}{\mathbb{R}}
\newcommand{\C}{\mathbb{C}}

\newcommand{\su}{\mathfrak{su}}

\newcommand{\SU}{{\rm SU}}

\newcommand{\SL}{\mathrm{SL}}

\DeclareMathOperator{\coker}{coker}

\DeclareMathOperator{\ind}{ind}
\DeclareMathOperator{\im}{im}

\DeclareMathOperator{\tr}{tr}

\renewcommand{\det}{\operatorname{det}}

\renewcommand{\epsilon}{\varepsilon}

\def\({\mathopen{}\left(}
\def\){\right)\mathclose{}}
\def\<{\mathopen{}\left<}
\def\>{\right>\mathclose{}}

\usepackage{multicol, color}

\definecolor{gold}{rgb}{0.85,.66,0}
\definecolor{cherry}{rgb}{0.9,.1,.2}
\definecolor{burgundy}{rgb}{0.8,.2,.2}
\definecolor{orangered}{rgb}{0.85,.3,0}
\definecolor{orange}{rgb}{0.85,.4,0}
\definecolor{olive}{rgb}{.45,.4,0}
\definecolor{lime}{rgb}{.6,.9,0}
\definecolor{green}{rgb}{.2,.7,0}
\definecolor{grey}{rgb}{.4,.4,.2}
\definecolor{brown}{rgb}{.4,.3,.1}

\def\makeautorefname#1#2{\AtBeginDocument{\expandafter\def\csname#1autorefname\endcsname{#2}}}

\newcommand{\mynewtheorem}[2]{
  \newaliascnt{#1}{equation}          
  \newtheorem{#1}[#1]{#2}
  \aliascntresetthe{#1}
  \makeautorefname{#1}{#2}
}
\mynewtheorem{theorem}{Theorem}
\mynewtheorem{prop}{Proposition}
\mynewtheorem{cor}{Corollary}
\mynewtheorem{construction}{Construction}
\mynewtheorem{lemma}{Lemma}
\mynewtheorem{conjecture}{Conjecture}

\numberwithin{substep}{step}
\makeautorefname{step}{Step}
\makeautorefname{substep}{Step}

\numberwithin{subcase}{case}
\makeautorefname{case}{Case}
\makeautorefname{subcase}{case}

\theoremstyle{remark}
\mynewtheorem{remark}{Remark}

\theoremstyle{definition}
\mynewtheorem{definition}{Definition}
\mynewtheorem{example}{Example}
\mynewtheorem{exercise}{Exercise}
\mynewtheorem{convention}{Convention}
\newtheorem*{convention*}{Convention}
\newtheorem*{conventions*}{Conventions}
\mynewtheorem{question}{Question}
\mynewtheorem{assumption}{Assumption}
        
\makeautorefname{chapter}{Chapter}
\makeautorefname{section}{Section}
\makeautorefname{subsection}{Section}
\makeautorefname{subsubsection}{Section}

\theoremstyle{introthm}
\newtheorem{introthm}{Theorem}
\newtheorem{introcor}[introthm]{Corollary}
\newtheorem{introconj}[introthm]{Conjecture}

\usepackage{amsrefs}
\usepackage{hyperref}

\newcommand\bbR{{\mathbb {R}}}

\newcommand\diamd{\blacklozenge}

\newcommand\chisu{\chi}
\def\F{\mathcal{F}}
\def\V{\mathcal{V}}
\def\CFloc{\mathit{CF}_{\operatorname{loc}}}
\def\Isloc{I^\sharp_{\operatorname{loc}}}
\def\Iloc{I_{\operatorname{loc}}}
\def\SI{\mathit{SI}}
\def\SIsloc{\SI^\sharp_{\operatorname{loc}}}
\def\colim{\operatorname{colim}}
\def\SU{\mathit{SU}}
\def\SL{\mathit{SL}}

\def\CS{\mathit{CS}}
\def\SF{\mathit{SF}}

\global\mdfdefinestyle{exampledefault}{% 
linecolor=black,backgroundcolor=gray!1,linewidth=1pt,leftmargin=0cm,rightmargin=0cm,topline=false,bottomline=false,skipabove=12pt}

\title{Nullhomotopic knots  are determined by their complements}

\author[Aliakbar Daemi]{Aliakbar Daemi}
\address{Department of Mathematics, Washington University in St. Louis\\
One Brookings Drive, St. Louis, MO 63130, USA}
\email{adaemi@wustl.edu}
\thanks{AD was partially supported by NSF grants DMS-1812033, DMS-2208181, DMS-2609113, NSF FRG Grant DMS1952762 and a Simons fellowship.} 

\author[Ciprian Manolescu]{Ciprian Manolescu}
\address{Department of Mathematics, Stanford University\\
450 Jane Stanford Way, Building 380, Stanford, CA 94305, USA}
\email{cm5@stanford.edu}
\thanks{CM was partially supported by a Simons Investigator Award, the Simons Collaboration Grant on New Structures in Low-Dimensional Topology, and the NSF AIMing grant DMS-2522743.}

\author[Nick Rozenblyum]{Nick Rozenblyum}
\address {Department of Mathematics, University of Toronto\\ 40 St. George Street, Room 6290, 
Toronto, Ontario M5S 2E4, Canada}
\email {nick.rozenblyum@utoronto.ca}
\thanks{NR was partially supported by an NSERC Discovery Grant (RGPIN-2025-0696).}

\begin{document}

\begin{abstract}
We prove that nullhomotopic knots in three-manifolds are determined by their complements.  This answers a Kirby Problem posed by Boileau. The argument uses techniques in instanton Floer homology and $SU(2)$-representation varieties.    
\end{abstract}
\maketitle

\hypersetup{linkcolor=black}
%\tableofcontents

%!TEX root = main.tex

\section{Introduction}\label{sec:intro}

The celebrated {\it Knot Complement Theorem} of Gordon and Luecke states that:

\begin{introthm}[\cite{GL}]\label{CKCP}
	If two knots $K$ and $K'$ in $S^3$ have homeomorphic exteriors, then there is a homeomorphism of $S^3$ taking $K$ to $K'$.
\end{introthm}

This result plays a fundamental role in low-dimensional topology, and consequently, it has become a major problem to extend this to other three-manifolds.  In fact, this is the focus of Problem 1.81 (D) in \cite{Kirby}, Conjecture 6.2 in \cite{GordonICM}, and Conjecture 1.1 in \cite{NiWu}:
\begin{introconj}\label{NKCP}
	Suppose $K$ and $K'$ are knots in a closed oriented 3-manifold $Y$ with orientation-preserving homeomorphic exteriors, such that these exteriors are irreducible and are not solid tori.  Then there exists an orientation-preserving homeomorphism of $Y$ mapping $K$ to $K'$.
\end{introconj}

There have been a few settings where the Knot Complement Theorem has been established, such as the lens spaces $L(p,q)$, when $p$ is square-free \cite{Gainullin}, and circle bundles over surfaces with genus at least 2 \cite{CremaschiYarmola}.  A general  
approach is to recast this as a Dehn surgery problem: if two knots $K$ and $K'$ in a three-manifold $Y$ have homeomorphic exteriors but there is no homeomorphism of $Y$ sending $K$ to $K'$, then that means that there is a non-trivial knot in $Y$ with a non-trivial surgery to $Y$.  To see this, we may assume that $K'$ is non-trivial. Applying the homeomorphism of the exterior sends the meridian of $K$ to a non-meridional slope on $K'$.  Doing that corresponding surgery on $K'$, then returns $Y$. Gordon and Luecke proved such non-trivial $S^3$ surgeries cannot exist using deep combinatorial techniques.  In other settings, various tools such as Heegaard Floer homology have been used for this.   

The main result of this paper is the following Dehn surgery problem.
 \begin{introthm}\label{thm:Boileau}
	Let $K$ be a non-trivial nullhomotopic knot in a closed, oriented three-manifold $Y$. A non-trivial surgery on 
	$K$ never gives a manifold homeomorphic to $Y$ in an orientation-preserving way.
\end{introthm}  

\begin{figure}
\begin{subfigure}{0.5\textwidth}
\labellist
\small\hair 2pt
\pinlabel {$5$} at 64 200
\pinlabel {$\textcolor{red}{K}$} at 200 195
\endlabellist
\centering
\includegraphics[scale=.5]{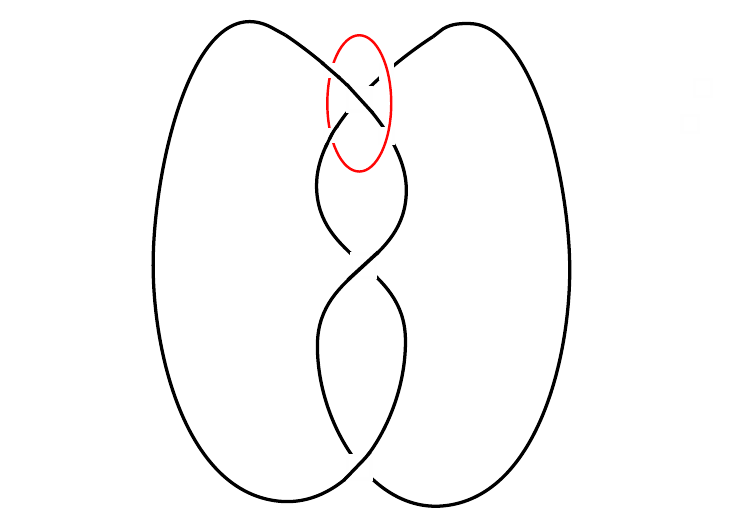}
\end{subfigure}
~
\begin{subfigure}{0.5\textwidth}
\labellist
\small\hair 2pt
\pinlabel {$5$} at 64 200
\pinlabel {$\textcolor{red}{K'}$} at 205 195
\endlabellist
\centering
\includegraphics[scale=.5]{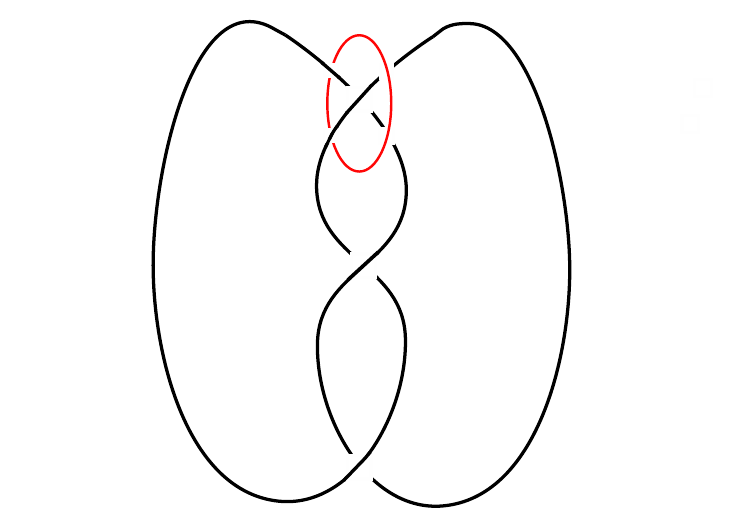}
\end{subfigure}
\caption{The left figure shows the knot $K$, a nullhomotopic knot in $L(5,1) = S^3_{5}(U)$.  (It is nullhomotopic since $K$ has linking number zero with the unknot, and in a lens space, nullhomologous knots are automatically nullhomotopic.)  Performing $1$-surgery on $K$ yields the figure on the right, where $K'$ is the core curve of the surgery.  The ambient manifold is $S^3_{5}(T_{2,3}) = -L(5,1)$, and this knot is again nullhomotopic.  It follows that $K$ and $K'$ have the same exterior, but the knots are not equivalent, being in different oriented three-manifolds.}    \label{lens}
\end{figure}

The statement of Theorem \ref{thm:Boileau} is essentially the same as \cite[Problem 1.80 (B)]{Kirby} proposed by Boileau. In the original statement of \cite[Problem 1.80 (B)]{Kirby}, existence of an orientation-preserving homeomorphism type is replaced with having the same simple homotopy type. It is shown in \cite[Theorem 1]{Turaev} that homeomorphism type is determined by the simple homotopy type in light of geometrization. On the other hand, it seems essential to add the orientation-preserving assumption on our homeomorphisms; it is observed in \cite{Morimoto} that there is a nullhomotopic hyperbolic genus one fibered knot $K$ in $L(5,1)$ with a surgery to $-L(5,1)$ (see Figure~\ref{lens}). Moreover, we drop the irreducibility assumption on $Y$ from \cite[Problem 1.80 (B)]{Kirby} in the statement of Theorem \ref{thm:Boileau}, since this assumption is not necessary. We also remark that, as in the Knot Complement Theorem, Theorem~\ref{thm:Boileau} readily implies the following {\em Nullhomotopic Knot Complement Theorem}:
 \begin{introthm}\label{thm:Boileau-2}
	If $K$ and $K'$ are knots in a closed, oriented 3-manifold $Y$ with orientation-preserving homeomorphic exteriors, and if $K$ is nullhomotopic, 
	then there exists an orientation-preserving homeomorphism of $Y$ mapping $K$ to $K'$.
\end{introthm}
\noindent This is the nullhomotopic case of Conjecture \ref{NKCP}, without the assumption that the exterior is irreducible and not a solid torus.

It is natural to ask about the following version of Theorem \ref{thm:Boileau-2}: If $K, K'$ are nullhomotopic knots in three-manifolds $Y, Y'$, then does having homeomorphic exteriors imply that the pairs $(Y, K)$ and $(Y',K')$ are {\em equivalent} (i.e., that there is a homeomorphism from $Y$ to $Y'$ taking $K$ to $K'$)?  This is Problem 1.80 (A) in \cite{Kirby} (with the additional assumption of irreducibility on $Y$ and $Y'$). The answer is negative in general.

\begin{example}
If one changes the framings in Figure~\ref{lens} from 5 to 7, then one obtains a nullhomotopic knot $K$ in $L(7,1)$ with a surgery to $L(7,4)$. This yields a pair of nullhomotopic knots with homeomorphic complements in non-homeomorphic lens spaces, independent of orientations. $\diamd$
\end{example}

However, if we impose the condition that one of the manifolds is not a lens space, then Problem 1.80 (A) in \cite{Kirby} has an affirmative answer, without even keeping track of orientations. 

\begin{introthm}
\label{thm:Boileau-3} 
Let $Y$ and $Y'$ be two closed, oriented three-manifolds. Let $K \subset Y$ and $K' \subset Y'$ be two nullhomotopic knots with homeomorphic complements. Assume that $Y$ is irreducible and not a lens space. Then the pairs $(Y, K)$ and $(Y', K')$ are equivalent.
\end{introthm}

\subsection{Sketch of the proof} 
\label{sec:sketch}
We tackle Boileau's conjecture (Theorem \ref{thm:Boileau}) using techniques from gauge theory.  Lackenby's work \cite[Theorem A.22]{Lackenby} implies Theorem \ref{thm:Boileau} holds in the case when the ambient three-manifold has $b_1 > 0$. Henceforth, we may restrict to the case of rational homology spheres. Our first result restricts the possible surgeries.  

\begin{introthm}\label{thm:|n|>1}
	Let $K$ be a non-trivial nullhomotopic knot in a (possibly reducible) oriented rational homology sphere $Y$.  If $m/n$-surgery on $K$ produces a three-manifold with the same fundamental group as $Y$, then $m/n = \pm 1$.         
\end{introthm}

Note that this theorem implies that Theorem \ref{thm:Boileau} holds if the surgery coefficient is not $\pm 1$ even if we drop the orientation-preserving assumption. In particular, the phenomenon provided by Figure \ref{lens} is specific to the case of surgery coefficient $\pm1$. As a consequence of Theorem \ref{thm:|n|>1}, we have the following result, which is a weaker version of Theorem~\ref{thm:Boileau-3} in general, but also applies to the case where $Y$ is reducible or a lens space.

\begin{introcor}\label{introcor}
Let $K$ be a non-trivial nullhomotopic knot in a rational homology sphere $Y$.  There is at most one nullhomotopic knot pair $(Y',K')$ inequivalent to $(Y,K)$ with homeomorphic complement.  
\end{introcor}

 In the statement of Corollary~\ref{introcor} and its proof, we do not assume that the involved homeomorphisms are orientation preserving.

\begin{proof}
 Suppose there exist two nullhomotopic knots in three-manifolds which are not equivalent to $K$, nor to each other, but they have homeomorphic complements.  These knots must arise as the core curves of non-trivial surgeries on $K$.  By Lemma~\ref{lem:nullhomotopic-core}, these surgered manifolds must have the same fundamental group as $Y$.  Theorem~\ref{thm:|n|>1} implies these knots arise as the core curves of $\pm 1$-surgery on $K$.  Denote these knots by $K_+$ and $K_-$ in $Y_+$ and $Y_-$ respectively, corresponding to the sign of the surgery coefficient.   We get that $-\frac{1}{2}$-surgery on $K_+$ in $Y_+$ is homeomorphic to $Y_-$.  However, since $K_+$ is nullhomotopic in $Y_+$, and $Y_+$ has the same fundamental group as $Y_-$, this contradicts Theorem~\ref{thm:|n|>1}.      
\end{proof}

A key ingredient in the proof of  Theorem \ref{thm:|n|>1} is the following result on detection of non-trivial nullhomotopic knots using instanton Floer homology.  
\begin{introthm} \label{non-triviality-I}
	Let $K$ be a non-trivial nullhomotopic knot in a rational homology sphere $Y$.  Then the instanton Floer homology group of the 0-surgery $Y_0(K)$ 
	is non-trivial.  
\end{introthm}

Since $Y_0(K)$ has a positive $b_1$, we need to work with instanton Floer homology for {\it admissible bundles} in this theorem. In particular, we need to specify a choice of an $SO(3)$-bundle on $Y_0(K)$. See Theorem \ref{non-triviality-I-detailed} for a more precise statement. The proof of Theorem \ref{non-triviality-I} combines several works including \cite{KMICM,DLVW:ribbon,HL:non-sep-S^3,Ni,DS:MCG-action}.

Returning to Theorem \ref{thm:|n|>1}, showing that $|m| = 1$ is trivial (see Lemma~\ref{lem:nullhomotopic-core}).  However, showing that $|n| = 1$ requires more work. Using Theorem \ref{non-triviality-I} and Lemma~\ref{lem:nullhomotopic-core}, it suffices to show that if for a nullhomotopic knot $K$, the core of $Y_{1/n}(K)$ with $|n|>1$ is nullhomotopic, then the instanton Floer homology of $Y_0(K)$ is trivial. We achieve this goal by picking an appropriate holonomy perturbation in the definition of the instanton Floer homology of $Y_0(K)$ such that the instanton Floer complex does not have any generators and hence it is trivial. In particular, our proof is modeled on the strategy used in \cite{PropertyPpillow} (see also \cite{Lin}).

Once Theorem \ref{thm:|n|>1} is established, to continue with the proof of Theorem~\ref{thm:Boileau} we can assume that the surgery coefficient is $-1$. (The case $n=1$ can be reduced to this case by switching the roles of $Y$ and $Y_{1}(K)$.)   

We proceed to outline the proof for $n=-1$. The argument we give later in the paper is not phrased in the same language as here, although it is essentially the same. 

For simplicity, in this introduction let us first also assume that $Y$ is an integer homology sphere.  Let $W$ denote the elementary two-handle cobordism from $Y$ to $Y_{-1}(K)$. Then Floer's surgery exact triangle \cite{floer:inst2, bd} asserts that the instanton Floer homology of $Y_0(K)$ is trivial if and only if the instanton Floer cobordism map $I(W):I(Y)\to I(Y_{-1}(K))$ is an isomorphism. 

The instanton Floer homology group of $Y$ is obtained by applying Morse homological methods to a Chern--Simons functional whose critical set can be identified with the $SU(2)$-character variety of $Y$ given as
\[
	\chi(Y):=\operatorname{Hom}(\pi_1(Y), SU(2))/{\sim}, 
\]
where $\sim$ denotes the conjugation action on the representation variety. To be a bit more specific, $I(Y)$ is the homology of a chain complex $(C(Y),d)$ with the chain group being $C(Y)$ generated by the non-trivial critical points of the Chern--Simons functional, possibly after a perturbation to guarantee the non-degeneracy of the critical points.  In particular, in the special case that we do not make any perturbation and the Chern--Simons functional is Morse, we can arrange so that $C(Y)$ is generated by $\chi^*(Y)$, the set of non-trivial elements of $\chi(Y)$. (Note that we might still need to perturb the Chern--Simons functional to define the differential $d$, which is given by the {\it downward gradient flow lines} of the Chern--Simons functional, or more precisely the solutions of the ASD equation on the cylinder $\mathbb{R}\times Y$.) 

More generally, the Chern-Simons functional can be used to provide a filtration on $(C(Y), d)$. In the associated graded, the differential only counts flow lines of small energy, and the resulting homology is called  {\em local instanton homology}, denoted $\Iloc(Y)$. We get a spectral sequence 
\begin{equation}\label{spec-seq}
  \Iloc(Y) \implies I(Y).
\end{equation}

A useful case to keep in mind is when the Chern--Simons functional is Morse-Bott. A three-manifold $Y$ for which this holds is called {\em $SU(2)$-non-degenerate}. In that case, $\chi^*(Y)$ is a smooth manifold (with components of possibly different dimensions), and $\Iloc(Y)$ is simply the homology of $\chi^*(Y)$. Indeed, after picking a small enough Morse function on $\chi^*(Y)$ then one can find a chain complex computing $I(Y)$ where, morally, the chain groups are generated by the critical points of the Morse function.

Given a negative definite cobordism $W:Y\to Y'$ between integer homology spheres $Y$ and $Y'$ with $b_1(W)=0$, there is a cobordism map $I(W):I(Y)\to I(Y')$. This homomorphism is defined in terms of the ASD equation, possibly after a perturbation. In fact, the cobordism $W$ induces a homomorphism of the spectral sequences \eqref{spec-seq} associated to $Y$ and $Y'$ such that the induced maps on the second and the final pages of the spectral sequence are respectively given by $\Iloc(W)$ and $I(W)$. In particular, if $\Iloc(W)$ is an isomorphism, then $I(W)$ is an isomorphism of instanton homology groups. This provides a useful tool to prove that a given cobordism map in instanton homology is an isomorphism because $\Iloc(W)$ is easier to analyze than $I(W)$.  A similar idea to study homomorphisms of Floer homology groups is discussed, for example, in \cite{DLVW:ribbon}.

It remains to show that when $W$ is the elementary cobordism $Y\to Y_{-1}(K)$, and $Y$, $Y_{-1}(K)$ are orientation-preserving homeomorphic, the map $\Iloc(W)$ is an isomorphism. The character variety of $W$ is in one-to-one correspondence with the moduli space of flat connections on $W$ and any such flat connection provides a solution of the ASD equation. Assuming that $Y$ and $Y':=Y_{-1}(K)$ are $SU(2)$-non-degenerate, one can show that the flat connections on $W$ provide regular solutions of the ASD equation. Furthermore, we obtain a correspondence of smooth manifolds
\begin{center}
  \begin{tikzpicture}
    \node (A) at (1,0) {$\chi^*(W)$};
    \node (B) at (-1,-2) {$\chi^*(Y)$};
    \node (C) at (3,-2) {$\chi^*(Y')$};
    \draw[->] (A) --  (B);
    \draw[->] (A) -- (C);
  \end{tikzpicture}
\end{center}
and the induced map of homology groups using the {\it pull-push} construction agrees with 
\[
  \Iloc(W):\Iloc(Y)\to \Iloc(Y').
\]
Using the assumption that $Y$, $Y_{-1}(K)$ are orientation-preserving homeomorphic, one may check that the maps in the correspondence are diffeomorphisms. In particular, the pull-push map is an isomorphism. This proves Theorem \ref{thm:Boileau} in the case of $SU(2)$-non-degenerate integer homology spheres. However, there exist integer homology spheres that are not $SU(2)$-non-degenerate; see \cite{BodenHeraldKirk}. In general, the character varieties may not be smooth, and the local instanton homologies do not have a simple description.  
 
The remedy is to deform the character variety, using a finite-dimensional family of perturbations. Each perturbation $h$ yields a space $\chi^*_h(Y)$, with $h=0$ giving back $\chi^*(Y)$. We can arrange so that the total space of the family, 
\begin{equation}\label{family-char-variety}{\bf X}^*(Y) = \bigcup_h\chi^*_h(Y),\end{equation} is smooth (and transversely cut-out) near $h=0$. We get a correspondence \begin{center}
  \begin{tikzpicture}
    \node (A) at (1,0) {${\bf X}^*(W)$};
    \node (B) at (-1,-2) {${\bf X}^*(Y)$};
    \node (C) at (3,-2) {${\bf X}^*(Y')$};
    \draw[->] (A) --  (B);
    \draw[->] (A) -- (C);
  \end{tikzpicture}
\end{center}
such that the maps in the correspondence are diffeomorphisms (near $h=0$). For generic $h$ near $0$, we have that $\chi^*_h(Y)$ is a finite set of points, and the push-pull construction puts it in bijection to $\chi^*_{h}(Y')$. The $h$-perturbed Chern-Simons functional gives a further filtration on the local instanton complex, such that the induced map on the associated graded is an isomorphism. From here we deduce that $\Iloc(W)$ and hence $I(W)$ are isomorphisms, completing the proof.

In order to include the more general case of rational homology spheres, we follow the same strategy but work with the framed version of instanton Floer homology from \cite{KMinst}, denoted $I^\sharp(Y)$, which is defined for any closed oriented 3-manifold $Y$. The role of $\chi(Y)$ is played by a framed character variety $\chi^\sharp(Y)$, which consists of two copies of the representation variety $R(Y)={\rm Hom}(\pi_1(Y),SU(2))$. We have a cobordism map
$$ I^\sharp(W) : I^\sharp(Y) \to I^{\sharp}(Y_{-1}(K))$$
and to show that it is an isomorphism, we prove it is so on the associated graded:
\begin{equation}
\label{eq:isloc}
 \Isloc(W) : \Isloc(Y) \to \Isloc(Y_{-1}(K)).
 \end{equation}
This can be done using the same method used for $\Iloc(W)$.

\subsection{Alternative proof} 
\label{sec:alternative}
We now sketch an alternative proof of Theorem~\ref{thm:Boileau}, which was actually discovered first. The beginning of the argument is the same, but to show that \eqref{eq:isloc} is an isomorphism we proceed differently. Rather than perturbing the Chern-Simons functional, we keep the (framed) character varieties unperturbed, and interpret them sheaf-theoretically. 

In fact, it is convenient to pass through the Atiyah-Floer correspondence \cite{DFL} and replace $I^\sharp$ with the framed symplectic instanton homology from \cite{WW}, denoted $SI^\sharp$. This is the Lagrangian Floer homology of two Lagrangians associated to a Heegaard splitting of the three-manifold. We need to prove that the analogue of \eqref{eq:isloc} for {\em local} symplectic instanton homology is an isomorphism:
\begin{equation}
\label{eq:sisloc}
 \SIsloc(W) : \SIsloc(Y) \to \SIsloc(Y_{-1}(K)).
 \end{equation}

To do this, we rely on our upcoming work \cite{DMR}, where we develop a general theory of sheaf quantization for local Lagrangian Floer cohomology. There, to well-behaved Lagrangians $L_0, L_1 \subset M$ we associate a local Floer cochain complex $\CFloc(L_0, L_1)$, whose generators are (roughly) the intersections of $L_0$ and a small Hamiltonian perturbation of $L_1$, and whose differential counts pseudo-holomorphic strips near $L_0 \cap L_1$.

We state here the main result from \cite{DMR}.

\begin{introthm}
\label{thm:sheaf}
 Let $M$ be a real analytic manifold, equipped with an analytic symplectic form. Let $L_0, L_1 \subset (M, \omega)$ be $C^p$-smooth (for some $p \geq 4$), subanalytic Lagrangian branes, and let $X = L_0 \cap L_1$ be their intersection.
 
To this data we associate a homotopy sheaf $\F_{L_0, L_1}$ on $X$ with the following properties:

\begin{enumerate}
\item[(a)] It only depends on local data near $X$; that is, if $(L_0', L'_1)$ is another Lagrangian pair that agrees with $(L_0, L_1)$ in a neighborhood of $X$, then $\F_{L_0, L_1}$ and $\F_{L'_0, L'_1}$ are quasi-isomorphic.

\item[(b)] Suppose $U \subset L_0 $ and $\Omega \subset M$ are open subsets, such that there is a subanalytic symplectomorphism $\psi: \Omega \to T^*U$ taking $\Omega \cap L_0$ to the zero section $U \subset T^*U$ and taking $\Omega \cap L_1$ to the graph of $df$ for some function $f: U \to \R$. Then the restriction of $\F_{L_0, L_1}$ to $U \cap X$ is isomorphic to a twist of the vanishing cycle sheaf $\V_f$ whose stalk at $x \in U$ is 
$$ (\V_f)_x = \underset{B}{\colim} \ C^*( \{y \mid f(y) \leq f(x)\} \cap B,  \{y \mid f(y) < f(x)\} \cap B; \rho),$$ 
where $B$ are the elements of any neighborhood basis of $x$ (ordered by inclusion). The twist $\rho$ is a graded rank-1 local system determined by the brane structure.

\item[(c)] If $X$ is compact, then the space of global sections $\F_{L_0, L_1}(X)$ is quasi-isomorphic to the local Floer cochain complex $\CFloc(L_0, L_1)$.
\end{enumerate}
\end{introthm}

An extension of Theorem~\ref{thm:sheaf} shows that we can interpret the map \eqref{eq:sisloc} as induced by a map between the vanishing cycle sheaves. To prove that a map of sheaves is a quasi-isomorphism, it suffices to check this on stalks. In our case, we use part (b) of Theorem~\ref{thm:sheaf} to reduce the condition on stalks to a question in singularity theory.

Theorem~\ref{thm:sheaf} has implications beyond Boileau's conjecture. The sheaf of vanishing cycles that we construct is a real-analytic analogue of the perverse sheaves inspired from Donaldson-Thomas theory that appear in complex and algebraic geometry \cite{BBDJS, BBBJstacks, Bussi}. In particular, we can apply our construction to the Lagrangian $G$-character varieties associated to a Heegaard splitting, for any reductive Lie group $G$. From here we obtain a $G$-Casson invariant of three-manifolds, and in fact a categorification of it (in the form of local Floer homology). This generalizes the classical $\SU(2)$ Casson invariant from \cite{Casson}, and also the $\SL(2, \C)$ Casson invariant and homology from \cite{AM}.

\subsection{Organization of the paper} In Section~\ref{sec:init}, we prove the relevant non-vanishing result, Theorem~\ref{non-triviality-I}.  We prove Theorem~\ref{thm:|n|>1} in Section~\ref{sec:|n|>1}.  Section~\ref{sec:|n|=1} is the main part of the paper; this is where we give a review of most of the relevant notions from instanton Floer homology and complete the proofs of Theorems~\ref{thm:Boileau}, \ref{thm:Boileau-2}, and \ref{thm:Boileau-3}.  In Section~\ref{sec:example}, we detail a technical point why the arguments fail if we do not require the homeomorphism between $Y_{-1}(K)$ and $Y$ to be orientation-preserving. 

\subsection{Statement of AI use} A previous version of this paper was written by the first author and Tye Lidman without any AI assistance; it contained a proof of the main theorem under the additional assumption that the three-manifolds are $SU(2)$-nondegenerate. For the past few years we have been writing (together with Tye Lidman) a complete proof of Boileau's conjecture along the lines of Subsection~\ref{sec:alternative}. In September 2026, we explored the capabilities of GPT-6 Astra in relation to Boileau’s conjecture, which led to the proof presented in this version of the paper (and outlined in Subsection~\ref{sec:sketch}). The key contribution from AI was the idea of using family versions of character varieties analogous to the one in \eqref{family-char-variety} to prove Lemma \ref{chi-h-chi'-h-smae-size}, which had previously been established under the additional assumption of $SU(2)$-non-degeneracy. Previously, the non-degeneracy assumption was also used, though less essentially, in the proofs of Lemma \ref{lem:index} and Proposition \ref{regular-almost-flat}. AI suggested adaptations of these proofs that would remove this assumption. Our revised proof of Lemma \ref{lem:index} follows essentially the suggested adaptation, while for Proposition \ref{regular-almost-flat}, we provide a different proof to avoid some of the technical difficulties in the AI-suggested argument.

We used AI for proofreading and for filling a gap in our proof of Theorem~\ref{thm:Boileau-3}. We also used it to produce Figure~\ref{immersed-disk} from a sketch. 

The arguments suggested by AI were checked and rewritten by the authors.  The article does not contain AI-generated text. 

We eventually plan to replace this version with one that contains our original proof, since we find it more enlightening. Tye Lidman will be added as a co-author to that version, but he declined to take part in the current write-up because of the AI involvement.

\subsection{Acknowledgments} We thank Marc Kegel, Tye Lidman and Yi Ni for helpful conversations.  
	
%!TEX root = main.tex

\section{Instanton homology and detection of the trivial knot}\label{sec:init}

The main goal of this section is to prove Theorem \ref{non-triviality-I} which gives a way to detect the unknot among null-homotopic knots using instanton homology. Before proving this theorem we review some of the basic properties of instanton homology. In particular, we make it more precise what version of instanton Floer homology is used in the statement of Theorem \ref{non-triviality-I}.

A pair of a 3-manifold $N$ and an embedded $1$-manifold $w$ in $N$, representing an element of $H_1(N;\Z)$, is called {\it admissible} if there is an embedded orientable surface $S$ in $N$ such that the algebraic intersection $w \cdot S$ is odd. This condition necessitates $b_1(N)>0$. Instanton Floer homology associates a relatively $\Z/8$-graded group to any admissible pair $(N,w)$ \cite{Floer}. This topological invariant of $(N,w)$ is defined as the Morse homology of a Chern--Simons functional (see subsection \ref{background-2}), and the critical points of this functional can be identified with
\begin{equation}\label{char-ad-pair}
  \chi(N,w):=\{\rho: \pi_1(N - w) \to SU(2) \mid \rho(\mu) = -1 \text{ for any meridian $\mu$ of $w$}\}/{\sim}.
\end{equation}
Here ${\sim}$ denotes the conjugation action. In particular, if $\chi(N,w)$ is empty, then $I(N,w) = 0$. It turns out that the instanton Floer homology only depends on the class of $w$ in $H_1(N; \Z/2)$. In the following, we will be working with the version of instanton Floer homology that is defined with coefficients in a field of characteristic zero (e.g. $\Bbb Q$).

\begin{example}\label{admisible-T3}
	A simple example of an admissible pair is given by $(T^3,S^1)$, where $S^1$ denotes one of the factors of $T^3=S^1\times S^1\times S^1$. 
	Then $I(T^3,S^1)$ has rank $2$, and the relative grading of the two generators for this instanton Floer homology group is $4$. Representatives for these two generators are given by the elements of 
	$\chi(T^3,S^1)$, which can be described in the following way. First note that any representation $\rho:\pi_1(T^2\setminus {\rm pt})\to SU(2)$ of the punctured torus that maps a small loop 	around the puncture to $-1$ is irreducible and can be conjugated so that the generators in the standard presentation of $\pi_1(T^2\setminus {\rm pt})$ are mapped to the elements $i,j\in SU(2)$. 
	The character variety $\chi(T^3,S^1)$ is equal to two copies of the character variety of the above punctured torus, each copy determined by the value of the representation at 
	the additional circle factor, which has to be a central element of $SU(2)$. $\diamd$
\end{example}

A key ingredient in the proof of Theorem \ref{non-triviality-I} is the following non-vanishing result. 

\begin{theorem}[\cite{KMsutures,DS:MCG-action}]\label{thm:km-non-vanishing}
	For an admissible pair $(N,w)$, the instanton homology group $I(N,w)$ is trivial if and only if there is an embedded sphere $S$ in $N$ such that $w\cdot S$ is odd. In particular, 
	$I(N,w)$ is non-trivial if $N$ is irreducible.
\end{theorem}

In the case that $N$ is irreducible, this result is due to Kronheimer and Mrowka \cite{PropertyPtriangle,KMsutures}. The latter proof, stated as \cite[Theorem 7.21]{KMsutures}, relies on the  machinery of instanton homology for sutured manifolds, developed in the same paper, and appears to require some modifications to both the statement and the proof. The general case of Theorem \ref{thm:km-non-vanishing} is discussed in \cite{DS:MCG-action}. This proof follows a plan proposed in \cite{DIS:rank3ins}, which uses a slight generalization of sutured instanton homology in \cite{KMsutures} and a connected sum theorem from \cite{Scaduto}. We also remark that an extension of Kronheimer and Mrowka’s result from the irreducible case to more general admissible pairs was also known to Tye Lidman and Juanita Pinz\'on-Caicedo.

Now we are ready to give a more precise version of Theorem \ref{non-triviality-I}.
\begin{theorem} \label{non-triviality-I-detailed}
	Let $K$ be a non-trivial nullhomotopic knot in a rational homology sphere $Y$.  Suppose $w_0$ is the embedded 1-manifold given by the core of surgery in $Y_0(K)$. Then, the instanton Floer homology groups 
	$I(Y_0(K),w_0)$ and $I(Y_0(K)\#T^3, w_0 \cup S^1)$ are both non-trivial.
\end{theorem}
\begin{proof}
This theorem follows from combining Theorem \ref{thm:km-non-vanishing} and the following result of Ni. 
\end{proof}
\begin{theorem}[\cite{Ni}]\label{thm:ni}
	If $K$ is a non-trivial nullhomotopic knot in a rational homology sphere, then 0-surgery cannot contain a non-separating 2-sphere.
\end{theorem}
In fact, Ni proves a more general result where the assumption that $Y$ is a rational homology sphere is replaced with the assumption that $Y\setminus K$ does not contain an $S^1\times S^2$ summand. The extension to the case $b_1(Y)>0$ follows from \cite{Lackenby,Gabai:fol-II}. In the case that $Y$ is a rational homology sphere, Ni uses some results about degree one maps of 3-manifolds and \cite{DLVW:ribbon} to conclude that if $Y_0(K)=Z\# S^1\times S^2$, then $\pi_1(Y)=\pi_1(Z)$. Theorem \ref{thm:ni} then follows from \cite{HL:non-sep-S^3}.

%!TEX root = main.tex

\section{Controlling the surgery coefficient}\label{sec:|n|>1}
In this section, we show that for a non-trivial nullhomotopic knot in a rational homology sphere, any surgery coefficient $m/n$ other than $+1$ or $-1$ has a homotopically essential dual knot. In the first subsection, we relate the fundamental group of the surgered manifold to the homotopy class of the dual knot. Next, we use the unknot detection result from the previous section to complete the proof of Theorem \ref{thm:|n|>1}.

\subsection{The case of $|m|\neq 1$}
We begin with the following topological lemma characterizing when surgery on a nullhomotopic knot can return a three-manifold with the same fundamental group.
\begin{lemma}\label{lem:nullhomotopic-core}
	Let $K$ be a nullhomotopic knot in $Y$.  Suppose that for the surgery coefficient $m/n$, we have $\pi_1(Y_{m/n}(K)) \cong \pi_1(Y)$.  Then, $m = \pm 1$ and the core of the Dehn surgery is also nullhomotopic.  The converse is also true; if the core of $m/n$-Dehn surgery on $K$ in $Y$ is also nullhomotopic, then $m = \pm 1$ and $\pi_1(Y_{m/n}(K)) \cong \pi_1(Y)$.  
\end{lemma}
\begin{proof}
	We begin with the forward direction.  Since $H_1(Y_{m/n}(K)) \cong H_1(Y) \oplus \Z/m$, we see that $m = \pm 1$.  To show that the core is nullhomotopic, we assume $m/n = 1/n$ with $n \geq 1$. The case of negative $n$ is similar.  
	Using the slam dunk move, we can view $1/n$-surgery on $K$ as $(+1)$-surgery along $K$ and a sequence of $(+2)$-surgeries along a chain of $n-1$ unknots.  
	Consider the corresponding 2-handle cobordism $W$ from $Y$ to $Y_{1/n}(K)$.  Since $K$ is nullhomotopic, we see that the map $\pi_1(Y)\to \pi_1(W)$ induced by inclusion is an isomorphism. We may similarly view $W$ 
	(with opposite orientation and flipped upside down) as attaching a 2-handle to $Y_{1/n}(K)$ along the core of the surgery and further 2-handles along a chain of unknots. 
	Therefore, the inclusion of $Y_{1/n}(K)$ in $W$ as the outgoing end induces a surjection from $\pi_1(Y_{1/n}(K)) \cong \pi_1(Y)$ to $\pi_1(W) \cong \pi_1(Y)$. 
	Since three-manifold groups are Hopfian, we see that the inclusion 
	from $Y_{1/n}(K)$ into $W$ induces an isomorphism on $\pi_1$. This means that the core must also be nullhomotopic.   

	To prove the converse, suppose that the core of surgery on $K$ is nullhomotopic.  Using a similar argument to above and identifying $Y_{m/n}(K)$ as integral surgery on $K$ and a chain of unknots, we have a 2-handle cobordism 
	$W : Y \to Y_{m/n}(K)$.  We have $\pi_1(W) \cong \pi_1(Y)$ since $K$ is 
	nullhomotopic.  By flipping $W$ upside down and using that the core is nullhomotopic, we get $\pi_1(Y_{m/n}(K)) \cong \pi_1(W)$.  Thus
	$\pi_1(Y_{m/n}(K))$ and $\pi_1(Y)$ are isomorphic, which also implies that $H_1(Y_{m/n}(K))\cong H_1(Y) \oplus \Z/m$ and $H_1(Y)$ are isomorphic.
	In particular, we have $m = \pm 1$.  
\end{proof}
The reference to the cobordism $W$ is not essential in the proof of the above lemma as the claims can be verified directly using Van Kampen’s theorem. However, similar cobordisms play a crucial role in Section~\ref{sec:|n|=1}, which is why we include the cobordism-based proof here.

\subsection{The case of $|n|\neq 1$}
In order to prove Theorem~\ref{thm:|n|>1}, we will need to study the character varieties of nullhomotopic knot complements using the {\it pillowcase}.  We give a quick recap of this theory here.  

\begin{figure}[ht!]
\centering
\includegraphics[scale=.3]{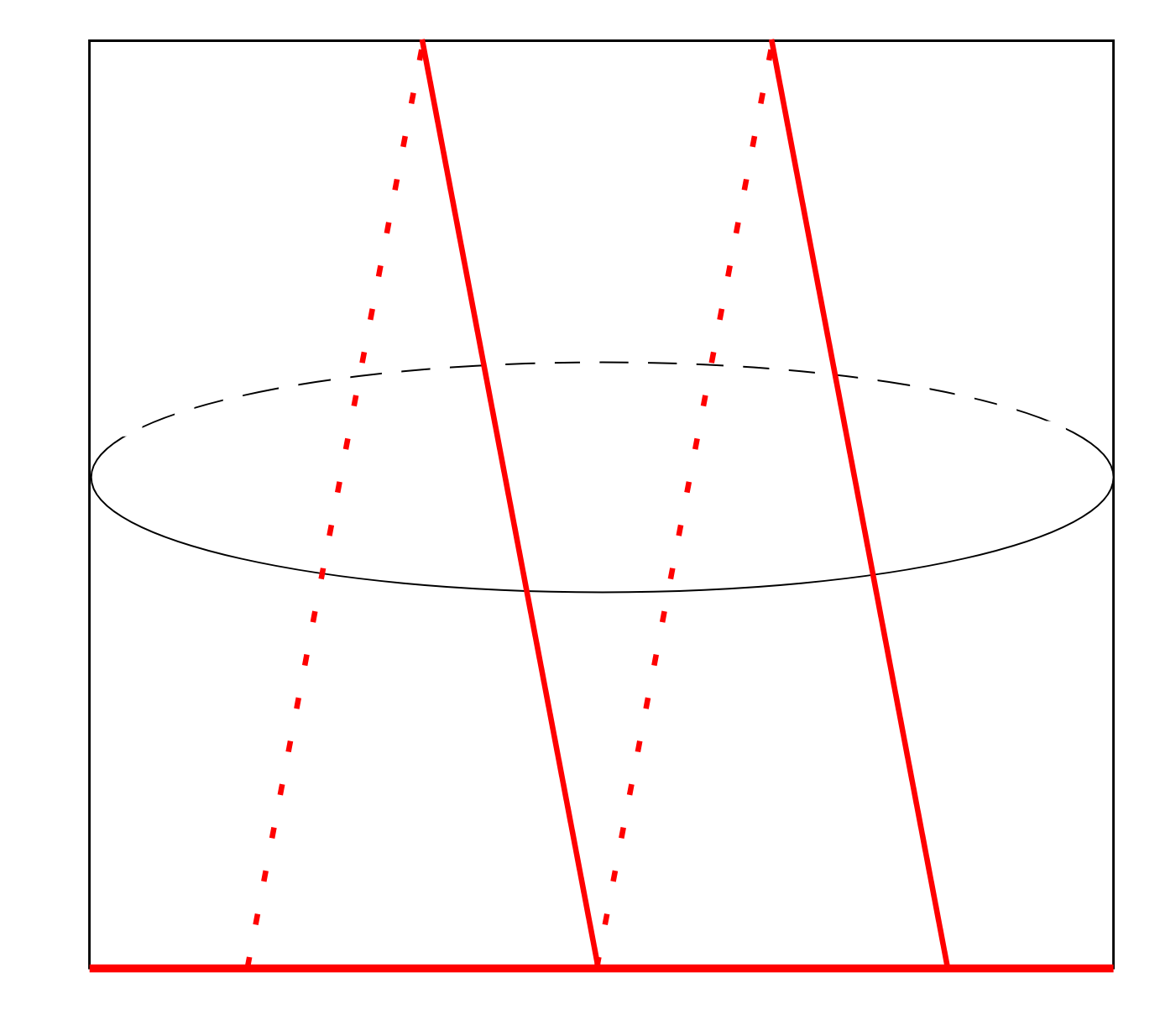}
\caption{The pillowcase is drawn in black.  In red is $\iota^*\chisu(E_{T_{2,3}})$, the image of the character variety for the right-handed trefoil under restriction.  The $x$-axis is given by the value of a character on $\mu$ and the $y$-axis is given by the value of a character on $\lambda$.}\label{fig:trefoil}
\end{figure}

Consider $\chisu(T^2)$, which is topologically the pillowcase orbifold, shown in Figure~\ref{fig:trefoil}.  Here, the coordinates are described as follows. Let $\mu,\lambda$ generate $\pi_1(T^2)$ and consider a representation $\rho : \mathbb{Z}^2 \to SU(2)$.  Since $\rho(\mu), \rho(\lambda)$ commute, we may conjugate them simultaneously to write 
\[
\rho(\mu) = \left[\begin{array}{cc} e^{i \alpha} & 0 \\ 0 & e^{- i\alpha} \end{array}\right], \hspace{1cm} 
\rho(\lambda) = \left[\begin{array}{cc}  e^{i \beta} & 0 \\ 0 & e^{- i\beta} \end{array}\right].
\] 
The pair $(\alpha, \beta) \in (\R/2\pi \Z)^2$ determines the conjugacy class of $\rho$, and the only other pair in $ (\R/2\pi \Z)^2$ representing the same conjugacy class is given by $(-\alpha,-\beta)$. Thus $\chisu(T^2)$ can be identified with the quotient of $(\R/2\pi \Z)^2$ with respect to the involution that maps $(\alpha, \beta)$ to $(-\alpha, -\beta)$. This involution has four fixed points $(0,0), (0,\pi), (\pi,0), (\pi,\pi)$, which give rise to  the  four orbifold points of the pillowcase $\chisu(T^2)$.
In our figures, we use the convention that the (conjugacy classes of) representations sending $\mu$ (resp. $\lambda$) to $1 \in SU(2)$ are those that comprise the ``left vertical seam'' (resp. ``bottom horizontal seam'') of the pillowcase in Figure~\ref{fig:trefoil}.

If $f : X \to Y$ is a map of spaces, then we obtain a pullback map $f^* : \chisu(Y) \to \chisu(X)$.  In particular, given a three-manifold with torus boundary, e.g. the exterior of a knot, we obtain a map from the character variety of the three-manifold to the pillowcase.  In general, this is {\em not} an embedding.  For notation, for a knot exterior $E_K$, we will write $\iota$ for the inclusion of $\partial E_K$ into $E_K$.  We will be interested in studying the shape of $\iota^*\chisu(E_K)$ inside the pillowcase.  Note that $\iota^*\chisu(E_K)$ is compact.  In the case of the exterior of a nullhomologous knot, any Seifert surface of $K$ determines a longitude $\lambda$ of $K$. In particular, we may use the meridian $\mu$ and the longitude $\lambda$ of $K$ to parametrize the pillowcase.  The image of $\iota^*\chisu(E_{T_{2,3}})$ is shown in Figure~\ref{fig:trefoil}.

We now illustrate why it is helpful to use the pillowcase picture—it provides a geometric intersection theory approach to computing character varieties. Suppose that $\iota^*\chi(E_K)$ intersects $\{ \alpha = 0\}$.  Then, this means there is an $SU(2)$ representation of $\pi_1(E_K)$ which sends $\mu$ to 1, in other words, there is an $SU(2)$ representation of $\pi_1(Y)$.  More generally, we obtain elements of $\chisu(Y_{p/q}(K))$ from intersections of $\iota^* \chi(E_K)$ with $\{(\alpha,\beta) \mid p\alpha+ q\beta= 0\}$.

We are now ready to see what the pillowcase can say about nullhomotopic knots.  
\begin{lemma}\label{lem:nullhomotopic}
Let $K$ be a nullhomotopic knot in a three-manifold $Y$.  Then, 
\begin{equation}\label{eq:nullhom-0}
\iota^*\chi(E_K) \cap \{\alpha = 0 \} = \{(0,0)\}
\end{equation}
and
\begin{equation}\label{eq:nullhom-pi}
\iota^*\chi(E_K) \cap \{\alpha = \pi\} = \{(\pi,0)\}.
\end{equation}
\end{lemma}
\begin{proof}
To prove \eqref{eq:nullhom-0}, suppose that $(0,x) \in \iota^*\chi(E_K)$.  Then, this means that there is an $SU(2)$-representation $\rho$ of $\pi_1(E_K)$ which sends $\mu$ to $I$ and $\lambda$ to $\operatorname{diag}(e^{ix}, e^{-ix})$.  Then $\rho$ induces a representation of $\pi_1(Y)$ sending $\lambda$ to $\operatorname{diag}(e^{ix}, e^{-ix})$.  Since $[\lambda] = [K] \in \pi_1(Y)$, and $K$ is nullhomotopic, $\rho(\lambda)=  I$.  Therefore, $x = 0$.  

Next, we turn to \eqref{eq:nullhom-pi}. Suppose that $(\pi, x) \in \iota^*\chi(E_K)$.  Since $K$ is nullhomotopic in $Y$, and $[\lambda] = [K] \in \pi_1(Y) = \pi_1(E_K) / \ll \mu \gg$, we have $\lambda \in \ll \mu \gg$ and hence we can write \[\lambda = g_1 \mu^{\epsilon_1} g_1^{-1} \cdots g_{k} \mu^{\epsilon_{k}} g_{k}^{-1},\]
where $\epsilon_i \in \{ \pm 1\}$.
Since $\lambda$ is nullhomologous in $E_K$ and $\mu$ is infinite order in $H_1(E_K)$, $k$ is an even integer.  Therefore, if $\rho(\mu) = -I$, then $\rho(\lambda) = I$.   
\end{proof}

\begin{remark}\label{rmk:core-surgery-pillowcase}
Suppose that the core $\widehat{K}$ of some $1/n$-surgery on $K$ is nullhomotopic.  Then by Lemma~\ref{lem:nullhomotopic}, 
\begin{equation}\label{eq:nullhom-surg}
\iota^*\chi(E_K) \cap \{\alpha + n \beta = 0 \} = \{(0,0)\}.
\end{equation}
This follows from \eqref{eq:nullhom-0} where we use a change of coordinates.  Indeed, $(\mu+n\lambda, \lambda)$ are the meridian and longitude for the core curve respectively, so in these coordinates for the pillowcase, $\iota^* \chi(E_{\widehat{K}}) \cap \{ \alpha = 0\} = \{(0,0)\}$ by Lemma~\ref{lem:nullhomotopic}.  Since $\iota^*\chi(E_K) = \iota^*\chi(E_{\widehat{K}})$, switching back to $(\mu,\lambda)$ coordinates gives \eqref{eq:nullhom-surg}. $\diamd$ 
\end{remark}

For 0-surgery on a nullhomologous knot, let $\mu_K$ be the dual knot of the surgery.  In that case, $\mu_K$ determines an admissible bundle on $Y_0(K)$.  Moreover, $\chi(Y_0(K), \mu_K)$ is empty if and only if $\iota^* \chi(E_K) \cap \{ \beta = \pi\}$ is empty.  Note that if $\chi(Y_0(K),\mu_K)$ is empty, then $I(Y_0(K),\mu_K) = 0$.  In fact, we have a more general vanishing result for $I(Y_0(K),\mu_K)$ that we can read off of the pillowcase.  

\begin{prop}\label{prop:nointersections}
Let $K$ be a nullhomologous knot in a rational homology sphere $Y$.  Let $\gamma$ be  an embedded path in the pillowcase from $(0,\pi)$ to $(\pi,\pi)$ which is homotopic to the horizontal path rel endpoints through paths whose interiors avoid the orbifold points.  If 
$\iota^*\chi(E_K) \cap \gamma = \emptyset$, then $I(Y_0(K),\mu_K) = 0$.  
\end{prop}
\begin{proof}
This result is well known.  For example, it is established for knots in integer homology spheres in \cite[Theorem 3.3]{LPZ} and for knots in rational homology spheres the argument generalizes without modification.  For other related results, see \cite{PropertyPpillow, Lin, Zentner}.  
\end{proof}
 
\begin{figure}[ht!]
\begin{center}
\includegraphics[scale=.4, trim = 0in 5.6in 0in 0in, clip]{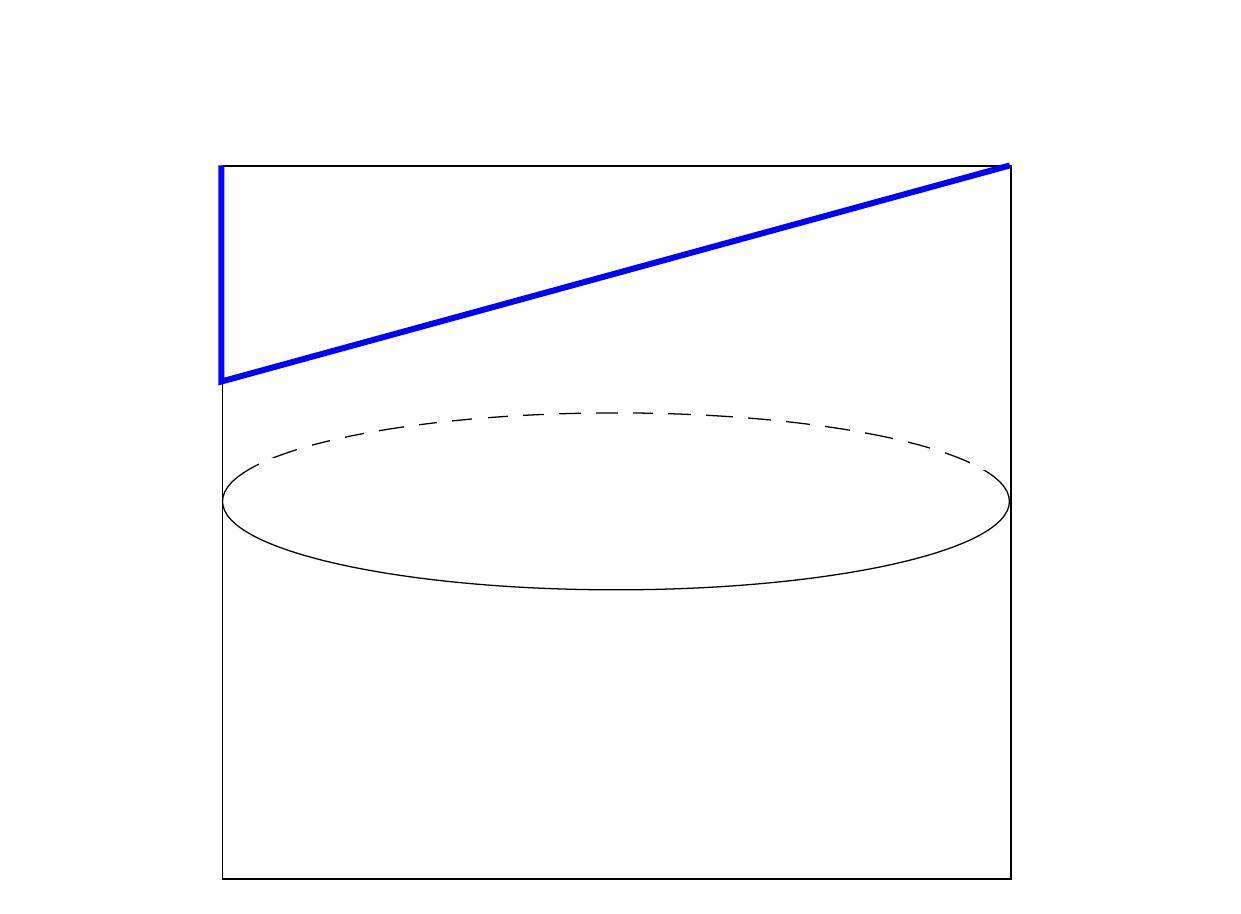}
\end{center}
\caption{The path $\gamma_3$ in the pillowcase.}\label{fig:holonomy-3}
\end{figure}

With the above work, we are now ready to prove Theorem~\ref{thm:|n|>1}, which we recall constrains when surgery on a nullhomotopic knot can return a three-manifold with the same fundamental group.  
\begin{proof}[Proof of Theorem~\ref{thm:|n|>1}]
Let $K$ be a non-trivial nullhomotopic knot in $Y$ such that $Y_{m/n}(K)$ has the same fundamental group as $Y$.  By Lemma~\ref{lem:nullhomotopic-core}, the surgery coefficient must be of the form $-1/n$ for some integer $n$.  By reversing orientation of $Y$ if necessary, we can assume that $n > 0$.  Suppose that $n \geq 2$ and we will obtain a contradiction.  Our goal is to find a path $\gamma$ from $(0,\pi)$ to $(\pi,\pi)$ in the pillowcase such that $\iota^*\chi(E_K) \cap \gamma$ is empty and then invoke Proposition~\ref{prop:nointersections}.

We define the path $\gamma_n$ from $(0,\pi)$ to $(\pi, \pi)$ as follows.  For odd $n$, this path is given by 
\[
\gamma_n = \{(0,\beta) \mid \frac{\pi(n-1)}{n} \leq \beta \leq \pi\} \cup \{(\alpha, \beta) \mid \alpha=n \beta-\pi(n-1) , \ \frac{\pi(n-1)}{n} \leq \beta \leq \pi \}
\]
For even $n$, define $\gamma_n$ by 
\[
\gamma_n = \{(\pi,\beta) \mid \pi \leq \beta \leq \frac{\pi(n+1)}{n}\} \cup \{(\alpha, \beta) \mid \alpha =n \beta-\pi n, \ \pi \leq \beta \leq \frac{\pi(n+1)}{n}  \}
\]
By construction, the path $\gamma_n$ is contained in the subset of the pillowcase represented by the union of the lines $\{\alpha = 0\}$, $\{\alpha = \pi\}$, and $ \{\alpha - n \beta = 0 \}$. Furthermore, for $n\geq 2$, it avoids the equivalence class of the point $(0,0)$ in the pillowcase.
The case $n = 3$ is shown in Figure~\ref{fig:holonomy-3}.  

By Lemma~\ref{lem:nullhomotopic} and Remark \ref{rmk:core-surgery-pillowcase}, $\iota^*\chi(E_K)$ is disjoint from $\gamma_n$ for $n\geq 2$.  Therefore, Proposition~\ref{prop:nointersections} implies that $I(Y_0(K),\mu_K) = 0$.  This contradicts Theorem~\ref{non-triviality-I-detailed}.
\end{proof}

It is natural to ask about why the argument above fails for $n = 1$.  The path $\gamma_{1}$ has an intersection with $\iota^*\chi(E_K)$ at $(0,0)$ coming, for example, from the trivial representation of $\pi_1(E_K)$ to $SU(2)$. One might still hope that there is a homotopy of $\gamma_1$ relative to its endpoints to get another curve $\gamma_1'$ that is disjoint from $\iota^*\chi(E_K)$. For example, if $\pi_1(Y)$ has no irreducible $SU(2)$ representations and $Y$ is an integer homology sphere (conjecturally $S^3$ is the only such 3-manifold), then there is $\epsilon>0$ such that the intersection of $\iota^*\chi(E_K)$ with the subset of the pillowcase that $\alpha\in [0,\epsilon)$ is contained in $\{\beta=0\}$, all induced by reducible representations on $E_K$. However, if $\pi_1(Y)$ has irreducible $SU(2)$ representations or $Y$ is not an integer homology sphere, then it is possible that $\iota^*\chi(E_K)$ separates $\{\alpha = 0\}$ and $\{\alpha -\beta = 0 \}$. For example, in a hypothetical case, $\iota^*\chi(E_K)$ could contain $\{\beta = 0\}$, $\{2 \alpha -\beta = 0\}$, but still misses $\{\alpha - \beta = 0 \}$, $\{\alpha = 0\}$ away from $(0,0)$. Then there would be no hope to find the desired $\gamma_1'$. Note that the proof of Theorem~\ref{thm:|n|>1} is not sensitive to whether a homeomorphism between $Y$ and $Y_1(K)$ is orientation preserving. In particular, the example in Figure \ref{lens} provides further evidence that adapting the above proof to $n=1$ would not be straightforward.

%!TEX root = main.tex
\section{Surgery slope $\pm 1$ and the surgery exact triangle}\label{sec:|n|=1}
\subsection{More background on instanton Floer homology of admissible pairs} \label{background-2}

Given a pair $(N,w)$ of a 3-manifold and a $1$-cycle, we can form the Hermitian rank $2$ bundle $E$ over $N$ determined uniquely up to isomorphism by the requirement that $c_1(E)$ is the Poincar\'e dual of $w$. Fix a $U(1)$ connection on $\Lambda^2 E$, and let $\mathcal A(E)$ be the space of connections on $E$ which induce the fixed connection on the determinant bundle $\Lambda^2 E$. The gauge group $\mathcal G(E)$ of the bundle automorphisms of $E$ with fiberwise determinant one acts on $\mathcal A(E)$. Let $\mathcal B(E)$ denote the quotient.

The Chern--Simons functional defines an $S^1$-valued function $\CS$ on $\mathcal B(E)$. After fixing a connection $B_0\in \mathcal A(E)$, define 
\[
  \CS(B):=\frac{1}{8\pi^2}\int_{[0,1]\times N} \tr\(F_{A}\wedge F_{A}\),
\]
where $A$ is a connection on the pullback of $E$ to $[0,1]\times N$ whose restrictions to $\{0\}\times N$ and $\{1\}\times N$ are respectively $B$ and $B_0$. This integral depends only on $B$ and $B_0$; applying an element of $\mathcal G(E)$ to $B$ changes $\CS$ by an integer and varying $B_0$ in the definition of $\CS$ changes it by an overall constant. In particular, we have the induced functional $\CS:\mathcal B(E)\to {\Bbb R}/{\Bbb Z}$, well-defined up to a constant shift. A connection $B$ on $E$ represents a critical point of $\CS$ if it is {\it projectively flat}, i.e., the adjoint connection associated to $B$ is flat. This characterization can be used to identify the critical set of $\CS$ with the character variety $\chi(N,w)$ in \eqref{char-ad-pair}. In the case that $(N,w)$ is admissible, all projectively flat connections on $E$ are irreducible.

The instanton Floer homology group $I(N,w)$ of an admissible pair can be regarded as the Morse homology of $\CS$. However, $\CS$ in general is not non-degenerate in the sense of Morse theory. In fact, the kernel of the Hessian of $\CS$ at a critical point $\alpha\in \chi(N,w)$ may be identified with the twisted cohomology group $H^1(N,{\rm ad}_\alpha)$, where ${\rm ad}_\alpha$ is the flat rank three vector bundle induced by $\alpha$ associated to the adjoint action of $U(2)$ on the Lie algebra $\frak {su}(2)$. To define $I(N,w)$ for an arbitrary admissible pair $(N,w)$, we need to perturb the Chern--Simons functional to obtain a non-degenerate functional.

A suitable family of such perturbations is provided by {\it cylinder functions}, and we review their definitions following \cite[Section 5.5]{donaldson-book}. First let $\gamma:S^1\times D^2 \to N$ be an embedded solid torus in $N$ and for any $z\in D^2$, let $\gamma_z$ be the loop given by the restriction of $\gamma$ to $S^1\times \{z\}$. For any connection $B$ representing an element of $\mathcal B(E)$, let ${\rm Hol}_z(B)$ be the holonomy of $B$ along $\gamma_z$.  Then the determinant of ${\rm Hol}_z(B)$ is determined by the fixed connection on $\Lambda^2 E$, and we fix a choice ${\rm d}_z(B)$ of a square root of $\det {\rm Hol}_z(B)$ that depends continuously with respect to $z$. Let ${\tau}_z(B)$ be the trace of ${\rm d}_z(B)^{-1}\cdot {\rm Hol}_z(B)\in SU(2)$. For a fixed smooth 2-form $\mu$ on $D^2$ with compact support and area one, we define 
\begin{equation}\label{special-cylinder}
  h(B):=\int_{D^2}\tau_z(B) \mu.
\end{equation}
Then $h$ induces a real valued map on $\mathcal B(E)$. More generally, a cylinder function is the composition of a smooth function $F:\Bbb R^n\to \Bbb R$ and $(h_1,\cdots,h_n):\mathcal B(E) \to \Bbb R^n$, where $h_i$ is given as in \eqref{special-cylinder}.

We may always arrange for a cylinder function $h$ such that the set of all critical points of the perturbed Chern--Simons function $\CS+\frac{1}{4\pi^2}h$, which is denoted by $\chi_h(N,w)$, consists of non-degenerate elements and $h$ is arbitrarily small \cite[Section 5.5.1]{donaldson-book}.  
For such a cylinder function $h$, let $C(N,w)$ be the vector space generated by $\chi_h(N,w)$ over a fixed coefficient field. After fixing a Riemannian metric on $N$, we may form the (formal) gradient of the perturbed Chern--Simons functional $\CS+\frac{1}{4\pi^2}h$, and the downward gradient flow lines for this functional can be identified with the solutions of the perturbed ASD equation over the cylinder $\Bbb R\times N$
\begin{equation}\label{perturbed-ASD-cyl}
  F_A^++(*_3\nabla_{A_t}h)^+=0,
\end{equation}
where we use the product metric on $\Bbb R\times N$ corresponding to the fixed metric on $N$. Here $A_t$ denotes the restriction of $A$ to $\{t\}\times N$, and $\nabla_{A_t}h$ defines a 1-form on $\{t\}\times N$. By applying Hodge star $*_3$ of the metric on $N$ and then putting together the resulting 2-forms on $\{t\}\times N$ for all values of $t$, we obtain a 2-form on $\Bbb R\times N$, and $(*_3\nabla_{A_t}h)^+$ is the self-dual part of this 2-form. 
Moreover, $F_A$ denotes the trace free part of the curvature of $A$ and $F^+_A$ is the self-dual part of this curvature. In particular, the left hand side of \eqref{perturbed-ASD-cyl} is a section of the vector bundle over $\Bbb R\times N$ with fiber $\Lambda^+\otimes \su(E)$, where $\su(E)$ denotes the bundle of trace free skew adjoint endomorphisms of $E$. For $\alpha,\beta\in \chi_h(N,w)$, we write $\breve M_h(N,w;\alpha,\beta)$ for the gauge equivalence classes of the solutions of \eqref{perturbed-ASD-cyl}, modulo the $\Bbb R$ action given by translation, which limit to perturbed projectively flat connections (gauge equivalent to) $\alpha$ as $t\to -\infty$ and $\beta$ as $t \to +\infty$. These moduli spaces for an appropriate choice of $h$ are cut out transversely, and $\breve M_h(N,w;\alpha,\beta)_d$ denotes the $d$-dimensional part of $\breve M_h(N,w;\alpha,\beta)$. After fixing a suitable orientation of these moduli spaces, for any $\alpha\in \chi_h(N,w)$, define 
\begin{equation}\label{differantial}
  d(\alpha)=\sum_{\beta\in \chi_h(N,w)} \#\breve M_h(N,w;\alpha,\beta)_0 \cdot \beta,
\end{equation}
and let $d:C(N,w)\to C(N,w)$ be the linear homomorphism determined by \eqref{differantial}. Then $d$ defines a differential on $C(N,w)$, and $I(N,w)$ is the homology of the chain complex $(C(N,w),d)$, which is an invariant of the homeomorphism type of $(N,w)$.

The instanton homology group $I(N,w)$ can be equipped with a relative $\Bbb Z/8$-grading, defined using the spectral flow of a family of self-adjoint elliptic operators parametrized by $\mathcal B(E)$. For any connection $B$ on $E$, consider the operator 
\[
  L_B^h:L^2_k(N,(\Lambda^0\oplus \Lambda^1)\otimes \su(E))\to L^2_{k-1}(N,(\Lambda^0\oplus \Lambda^1)\otimes \su(E))
\]
defined as 
\[
  L^h_B:=
  \left[
  \begin{array}{cc}
  	0&-d_B^*\\
	-d_B&*d_B+{\rm Hess}_{B}h
  \end{array}
  \right],
\]
where ${\rm Hess}_{B}h$ is the Hessian of $h$ at $B$. (We drop $h$ from our notation for $L^h_B$, if it is trivial.) This operator may be regarded as an extended Hessian of the perturbed Chern--Simons functional $\CS+\frac{1}{4\pi^2}h$. To see this, consider the following equation for $(B,\xi)\in \mathcal A(E)\times L^2_k(N,\su(E))$:
\begin{equation}\label{ext-pert-CS}
  *F_B+\nabla_Bh-d_B\xi=0
\end{equation}
This equation is invariant with respect to the action of the gauge group $\mathcal G(E)$ induced by the standard actions of the gauge group on $\mathcal A(E)$ and the sections of $\su(E)$. Furthermore, the term $*F_B+\nabla_Bh$, which is the gradient of $\CS+\frac{1}{4\pi^2}h$, is orthogonal to any term of the form $d_B\xi$ because of the gauge equivariance of  $\CS+\frac{1}{4\pi^2}h$. Therefore, any solution of \eqref{ext-pert-CS} satisfies $*F_B+\nabla_Bh=d_B\xi=0$. In particular, $[B]\in \chi_h(N,w)$, and then the irreducibility of $B$ implies that $\xi$ is trivial. The second row of $L^h_B$ gives a linearization of \eqref{ext-pert-CS} and the first row corresponds to the Coulomb gauge slice for the action of the gauge group on $ \mathcal A(E)\times L^2_k(N,\su(E))$. Furthermore, for $[B]\in \chi(N,w)$, the kernel of $L_B$ may be identified with $H^1(N,{\rm ad}_{B})$.

For any pair of connections $B,B'\in \mathcal A(E)$, the spectral flow of the family of the operators $\{L^h_{B(t)}\}$ associated to a path of connections $\{B(t)\}$ from $B$ to $B'$ determines an integer $\SF_h(B,B')$, which depends only on $B$, $B'$, $h$ and not on the specific path between $B$ and $B'$ because $\mathcal A(E)$ is contractible. In the case that $L^h_B$ and $L^h_{B'}$ have trivial kernel, then $\SF_h(B,B')$ is equal to the index of the ASD operator 
\[
  \mathcal D^h_A: L^2_k(\mathbb{R} \times N,\Lambda^1\otimes \su(E)) \to L^2_{k-1}(\mathbb{R} \times N,(\Lambda^0\oplus \Lambda^+)\otimes \su(E))
\]
defined as
\begin{equation}\label{DAh-cyl}
  \mathcal D^h_A(\zeta):=(-d_A^*\zeta,d_A^+\zeta+(*_3{\rm Hess}_{A_t}h(\zeta))^+),
\end{equation}
where $A$ is a connection on $\Bbb R\times N$ asymptotic to $B$ and $B'$ as $t\to -\infty$ and $\infty$. Here $*_3{\rm Hess}_{A_t}h(\zeta)$ is defined in a similar manner as $*_3\nabla_{A_t}h$: for any $t\in \Bbb R$, we take the restriction of $\zeta$ to $\{t\}\times N$, apply the operator $*_3{\rm Hess}_{A_t}h$ to obtain a 2-form on $N$, and then put together these 2-forms for all values of $t$ to produce a 2-form on $\Bbb R\times N$. The first component of \eqref{DAh-cyl} is given by the Coulomb gauge condition and the second component is obtained by linearizing \eqref{perturbed-ASD-cyl}.

For $\alpha, \alpha'\in \mathcal B(E)$, we may also consider $\SF_h(B,B')$ for representatives $B$ of $\alpha$ and $B'$ of $\alpha'$. Since $\pi_1(\mathcal B(E))=\Bbb Z$ and the spectral flow around a loop representing a generator of $\pi_1(\mathcal B(E))$ is $8$, the mod $8$ value of $\SF_h(B,B')$ is independent of the choice of the representatives $B$, $B'$ and is denoted by $\SF_h(\alpha,\alpha')$. Let $\theta_0\in \mathcal B(E)$ be the class represented by $B_0$, and define the grading of $\alpha\in \chi_h(N,w)$ as $\SF_h(\alpha,\theta_0)$. This grading on $C(N,w)$ is independent of $\theta_0$ up to an overall shift by a constant integer. Furthermore, the differential has degree $-1$ with respect to this grading because of the relation between the relative grading and the index of the ASD operator \eqref{DAh-cyl}.

Instanton Floer homology is functorial with respect to cobordisms. A cobordism $(W,c):(N,w)\to (N',w')$ of admissible pairs consists of a 4-dimensional cobordism $W$ of 3-manifolds and a relative 2-cycle $c$ in $W$, which can be used to define a Hermitian rank $2$ bundle on $W$ analogous to the 3-dimensional case. We fix a Riemannian metric on $W$, which is the product metric in a tubular neighborhood of its boundary corresponding to fixed metrics on $N$ and $N'$. Then we add two half-cylinders $(-\infty,-1]\times N$ and $[1,\infty)\times N'$ with product metrics to $W$ to obtain a Riemannian manifold with cylindrical ends. With a slight abuse of notation, we still write $W$ for this resulting non-compact Riemannian manifold. Suppose $h$, $h'$ are cylinder functions for $(N,w)$ and $(N',w')$ such that the corresponding perturbed Chern--Simons functionals have non-degenerate critical points. Then we can form the following perturbed ASD equation on $W$ 
\begin{equation}\label{perturbed-ASD-sec}
  F_A^++\rho_-(t)(*_3\nabla_{A_t}h)^++\rho_+(t)(*_3\nabla_{A_t'}h')^++\omega(A)=0,
\end{equation}
where $A_t$ for $t\in (-\infty,-1]$ and $A_t'$ for $t\in [1,\infty)$ are respectively restrictions of $A$ to $\{t\}\times N$, $\{t\}\times N'$, the terms $(*_3\nabla_{A_t}h)^+$ and $(*_3\nabla_{A_t'}h')^+$ are defined as in \eqref{perturbed-ASD-cyl}, the smooth function $\rho_-:(-\infty,-1]\to \Bbb R$ is  supported on $(-\infty,-2]$ and is equal to $1$ on $(-\infty,-3]$, and the smooth function $\rho_+:[1,\infty)\to \Bbb R$ is supported on $[2,\infty)$ and is equal to $1$ on $[3,\infty)$. Moreover, $\omega(A)$ is a secondary holonomy perturbation term which depends only on the restriction of $A$ to a fixed compact subspace of $W$ and it is added to the equation to give us enough flexibility to obtain the desired transversality results.

Let $\pi=(\omega,h,h')$ and write $M_\pi(W,c;\alpha,\alpha')$ for the moduli space of solutions to \eqref{perturbed-ASD-sec} that are asymptotic to $\alpha\in \chi_{h}(N,w)$ and $\alpha'\in \chi_{h'}(N',w')$ on the ends. Our assumption implies that $M_\pi(W,c;\alpha,\alpha')$ is cut out by a (non-linear) Fredholm equation and we write $M_\pi(W,c;\alpha,\alpha')_d$ for the components with expected dimension $d$. For any connection $A$ representing an element of $M_\pi(W,c;\alpha,\alpha')$, we can form the ASD operator $\mathcal D_A^\pi$ analogous to \eqref{DAh-cyl}, which is given by the pair of $-d_A^*$ and the linearization of \eqref{perturbed-ASD-sec} at $A$. By picking $\omega$ appropriately, we may assume that for any such $A$, the ASD operator $\mathcal D_A^\pi$ is surjective and hence $M_\pi(W,c;\alpha,\alpha')_d$ is a smooth $d$-dimensional manifold. In fact, this manifold is orientable and after a suitable orientation, we define 
\[
  C(W,c):C(N,w) \to C(N',w')
\]
as 
\[
  C(W,c)(\alpha)=\sum_{\alpha'\in \chi_{h'}(N',w')} \#M_\pi(W,c;\alpha,\alpha')_0 \cdot \alpha'.
\]
This defines a chain map respecting the relative gradings on $C(N,w)$ and $C(N',w')$. The induced map $I(W,c):I(N,w)\to I(N',w')$ at the level of homology groups does not depend on the choice of the metric and the perturbation term $\pi$ and is well-defined up to a sign.

Instanton Floer homology for admissible pairs can be used to define an invariant of arbitrary 3-manifolds.  Let $(T^3, S^1)$ be the admissible pair provided by Example \ref{admisible-T3}. The connected sum of this pair and an arbitrary $Y$ gives the admissible pair $Y^\sharp:=(Y \# T^3, S^1)$. Then the {\it framed instanton Floer homology} of $Y$ is defined as $I^\sharp(Y):=I(Y \# T^3, S^1)$; see \cite{KMinst}. Even though it is not relevant for this paper, we remark that we need to fix a basepoint of $Y$ together with a framing of its tangent space to keep track of where the connected sum is performed so that $I^\sharp(Y)$ is an invariant of $Y$ up to a canonical isomorphism. 

The critical points of the Chern--Simons functional for the admissible pair $Y^\sharp$ can be characterized in terms of the representation variety of $Y$. More specifically, we have 
\begin{equation}\label{char-variety-sharp=rep-variety}
  \chi(Y^\sharp)=R(Y)\sqcup R(Y),
\end{equation}
where $R(Y)={\rm Hom}(\pi_1(Y),SU(2))$. In particular, $\chi(Y^\sharp)$ admits the structure of a real algebraic variety. The identification in \eqref{char-variety-sharp=rep-variety} follows from the observation, recalled in Example \ref{admisible-T3}, that $\chi(T^3, S^1)$ consists of two irreducible representations. Since these two elements are irreducible and non-degenerate, a straightforward Mayer--Vietoris argument applied to $Y \# T^3$ gives the following exact sequence for any $\alpha\in \chi(Y)$ and an element $\alpha^\sharp\in \chi(Y^\sharp)$ obtained by gluing $\alpha$ to one of the two elements of $\chi(T^3, S^1)$:
\begin{equation}\label{exactsequence-Zariski-tangent-space}
	0\xrightarrow{} H^0(Y,{\rm ad}_\alpha)\xrightarrow{} \mathfrak{su}(2)\xrightarrow{} H^1(Y\# T^3,{\rm ad}_{\alpha^\sharp})\xrightarrow{} H^1(Y,{\rm ad}_{\alpha})\xrightarrow{} 0.
\end{equation}
Here the third term is the cohomology group of the splitting sphere $S^2$ in $Y \# T^3$ with trivial coefficient system $\mathfrak{su}(2)$.
To give a geometric explanation for \eqref{exactsequence-Zariski-tangent-space}, note that gluing $\alpha\in \chi(Y)$ to the elements of $\chi(T^3, S^1)$ contributes two submanifolds of dimension $3-\dim H^0(Y; {\rm ad}_{\alpha})$, which are either points, 2-dimensional spheres or 3-dimensional real projective spaces respectively when $\dim H^0(Y; {\rm ad}_{\alpha})$ is $3$, $1$ or $0$. Now deforming $\alpha$ within $\chi(Y)$ gives the remaining degrees of freedom in $H^1(Y\# T^3,{\rm ad}_{\alpha^\sharp})$ to deform $\alpha^\sharp$ corresponding to the elements of $H^1(Y,{\rm ad}_{\alpha})$.  

Similarly, if $W:Y\to Y'$ is a cobordism of based 3-manifolds together with a framed path connecting the basepoints, we can form the path sum of the product cobordism $I\times (T^3, S^1)$ and $W$ to form $W^\sharp:Y^\sharp\to {Y'}^\sharp$. Here the sum is performed in a neighborhood of the path connecting the basepoints. The cobordism map associated to $W^\sharp$ defines the cobordism map $I^\sharp(W):I^\sharp(Y) \to I^\sharp(Y')$. 

\subsection{Gradient-like perturbations of the ASD equation and perturbed flat instantons}\label{ASD-pert-gr-type}

There is an energy filtration on instanton Floer complexes of  admissible pairs that plays a crucial role in our proof of Theorem \ref{thm:Boileau} where the surgery coefficient is $\pm 1$. However, the general scheme of perturbing the ASD equation as in \eqref{perturbed-ASD-sec} does not interact nicely with this filtration. In this subsection, we introduce {\it gradient-like perturbations} of the ASD equation so that the corresponding cobordism maps behave better with respect to the energy filtration. 

Suppose $(N,w)$ and $(N',w')$ are admissible pairs and $(W,c):(N,w) \to (N',w')$ is a cobordism of admissible pairs. We fix a Riemannian metric on $W$ with cylindrical ends. For a (not necessarily closed) Riemannian 3-manifold $U$, we say that $W$ contains a cylinder modeled on $U$ if there is an isometric embedding $\iota: \Bbb R\times U\to W$ that is compatible with the cylindrical ends of $W$. The last condition asserts that there are isometric embeddings $i:U\to N$ and $i':U\to N'$ such that $\iota$ maps $(t,x)\in (-\infty,-1)\times U$ (resp. $(t,x)\in(1,\infty)\times U$) to $(t,i(x))$ (resp. $(t,i'(x))$). In the following, we will regard $\Bbb R\times U$ as a subspace of $W$ and we drop $\iota$ from our discussion. 

\begin{example} \label{example-cylinder-surgery-cob}
	Fix an admissible pair $(N,w)$ and a knot $K\subset N$ away from $w$. Let $N'=N_{-1}(K)$ and $W:N\to N'$ be the 2-handle cobordism. Since $w$ is away from $K$,
	it induces $w'\subset N'$ and a 2-cycle $c\subset W$. Even though it is not important for this example, we assume that $(N',w')$ is admissible so that it fits into the above setup. 
	For any codimension zero submanifold $U$ of the knot exterior, there is a natural embedding 
	$\iota: \Bbb R\times U\to W$, and we can arrange for a metric on $W$ with cylindrical ends which contains a cylinder modeled on $U$ corresponding to $\iota$. $\diamd$
\end{example}

Let $(W,c)$ contain a cylinder modeled on $U$ and $h:\mathcal A(U)\to \Bbb R$ be a cylinder function. The gradient-like perturbation of the ASD equation on $(W,c)$ associated to $h$ is defined as
\begin{equation}\label{perturbed-ASD}
  F_A^++(*_3\nabla_{A_t}h)^+=0,
\end{equation}
where $A$ is a connection on the bundle associated to $(W,c)$, $A_t$ denotes the restriction of $A$ to $\{t\}\times U$ and the second term in \eqref{perturbed-ASD} is defined in the same way as in \eqref{perturbed-ASD-cyl}.  Then \eqref{perturbed-ASD} is equivariant with respect to the action of the gauge group.

Define the {\it $h$-perturbed topological energy} for a connection $A$ on $(W,c)$ as 
\begin{equation}\label{top-energy-h}
	\mathcal E_h(A):=\frac{1}{8\pi^2}\int_W\tr\left((F_A+*_3\nabla_{A_t}h)\wedge (F_A+*_3\nabla_{A_t}h)\right).
\end{equation}
A similar energy function is defined in \cite[Section 3]{DFL}. Assuming the connection $A$ is decaying fast enough and is asymptotic to connections $\alpha$, $\alpha'$ on the ends, then one can easily check that 
\begin{equation}\label{top-energy-h-rewrite}
	\mathcal E_h(A)=\frac{1}{4\pi^2}(h(\alpha)-h(\alpha'))+\frac{1}{8\pi^2}\int_W\tr\left(F_A\wedge F_A\right).
\end{equation}
Here $h(\alpha)$ (resp. $h(\alpha')$) is defined by applying $h$ to the restriction of $\alpha$ to $U\subset N$ (resp. $U\subset N'$). The second term on the right hand side of \eqref{top-energy-h-rewrite} corresponds to what is commonly referred to as topological energy. We may also rewrite the $h$-perturbed topological energy as
\[
  \mathcal E_h(A)=\frac{1}{8\pi^2}\left (\left \vert\!\left \vert F_A+*_3\nabla_{A_t}h\right \vert\!\right \vert_{L^2}^2-
  2\left \vert\!\left \vert (F_A+*_3\nabla_{A_t}h)^+\right \vert\!\right \vert^2_{L^2}\right).
\]
This immediately gives the following lemma.
\begin{lemma}\label{pos-top-energy-sol}
	For any solution $A$ of \eqref{perturbed-ASD}, the $h$-perturbed topological energy $\mathcal E_h(A)$ is non-negative. Furthermore, for any such $A$, we
	have $\mathcal E_h(A)=0$ if and only if it satisfies the equation 
	\begin{equation}\label{pert-flat}
	  F_A+*_3\nabla_{A_t}h=0.
	\end{equation}   
\end{lemma}

We call any solution $A$ of \eqref{pert-flat} an $h$-perturbed flat connection. (A more accurate name is $h$-perturbed projectively flat connection, but we use the term $h$-perturbed flat connection instead for brevity.) The space of gauge equivalence classes of $h$-perturbed flat connections is the counterpart of the space of {\it constant pairs} in \cite{DFL}. An $h$-perturbed flat connection on the cobordism $(W,c):(N,w) \to (N',w')$ can be characterized in the following way. After applying a gauge transformation to an $h$-perturbed flat connection $A$ on $(W,c)$, we may assume that $A$ is in temporal gauge on the cylindrical end $(-\infty,-1]\times N$ and $\Bbb R\times U$. Then there is a connection $B$ on $(N,w)$ that satisfies the equation 
\begin{equation}\label{per-CS-crit}
  F_B+*_3\nabla_B h=0
\end{equation}
such that the restriction of $A$ to $Z=(-\infty,-1]\times N\cup \Bbb R\times U$ is the pullback of $B$ with respect to the projection map from $Z$ to $N$. (Note that \eqref{per-CS-crit} is equivalent to saying that $B$ represents an element of $\chi_h(N,w)$.) On the complement of $Z$, the connection $A$ is flat and one can identify $A$ in terms of its monodromies  along homotopy classes of loops.  

\begin{example}\label{elem-cob-almost-flat}
	Let $(N,w)$ be an admissible pair and $K\subset N$ be a knot away from $w$. 
	Let $(W,c):(N,w) \to (N',w')$ be the cobordism from Example \ref{example-cylinder-surgery-cob} with a cylinder modeled on $U$, a codimension zero submanifold of the knot 
	exterior. For any cylinder function $h$ on $U$, the above characterization of $h$-perturbed flat connections allows us to identify the moduli of $h$-perturbed flat connections on 
	$(W,c)$ with the subspace of $\chi_h(N,w)$ represented by connections that have trivial monodromy along $K$. $\diamd$
\end{example}

Next, we consider the local behavior of moduli spaces of solutions to \eqref{perturbed-ASD} around an $h$-perturbed flat connection. In this case, the ASD operator for a connection $A$ has the form 
\[
  \mathcal D^h_A: L^2_k(W,\Lambda^1\otimes \su(E)) \to L^2_{k-1}(W,(\Lambda^0\oplus \Lambda^+)\otimes \su(E))
\]
defined as
\begin{equation}\label{DAh}
  \mathcal D^h_A(\zeta):=(-d_A^*\zeta,d_A^+\zeta+(*_3{\rm Hess}_{A_t}h(\zeta))^+),
\end{equation}
where the term $(*_3{\rm Hess}_{A_t}h(\zeta))^+$ in \eqref{DAh} is defined in the same way  as in \eqref{DAh-cyl}.
The following lemma provides a description of the kernel of the ASD operator for $h$-perturbed flat connections and is a counterpart of the results of \cite[Section 4.3]{DFL}.

\begin{lemma}\label{char-ker-coker}
	Suppose $A$ is an $h$-perturbed flat connection on $(W,c)$. Then the kernel of the ASD operator $\mathcal D^h_A$ is given by
	\begin{equation}\label{ker-lin-const}
		\ker(\mathcal D^h_A)=\{\zeta \in L^2_k(W,\Lambda^1\otimes \su(E))\mid d_A^h\zeta=0,\, d_A^*\zeta=0\},
	\end{equation}
	where $d_A^h\zeta$ is the element of $L^2_{k-1}(W,\Lambda^2\otimes \su(E))$ given as $d_A\zeta+*_3{\rm Hess}_{A_t}h(\zeta)$.
\end{lemma}
\begin{proof}
	The claim about the kernel of $\mathcal D^h_A$ follows from the identity 
	\begin{equation}\label{h1-vanish}
	  \int_{W}\tr(d_A^h\zeta\wedge d_A^h\zeta)=0
	\end{equation}
	for any $\zeta\in L^2_1(W,\Lambda^1\otimes \su(E))$. In the simpler case that $h=0$ (in particular, $A$ is flat), 
	we have
	\begin{align}
		\int_{W}\tr(d_A\zeta\wedge d_A\zeta)&=\int_{W}d\tr(\zeta\wedge d_A\zeta)+\int_{W}\tr(\zeta\wedge d_A(d_A\zeta))\nonumber\\
		&=\int_{W}d\tr(\zeta\wedge d_A\zeta)+\int_{W}\tr(\zeta\wedge F_A\wedge \zeta),\label{linear-CW}
	\end{align}	
	where the first term in \eqref{linear-CW} vanishes by Stokes theorem and our assumption on the asymptotic 
	behavior of $\zeta$, and the second term in \eqref{linear-CW} vanishes because $A$ is flat.
	The general case of \eqref{h1-vanish} can be proved by analyzing the contributions from $\Bbb R\times U$
	and its complement to the integral in \eqref{h1-vanish} separately using again Stokes theorem 
	and then observing that the contributions cancel each other out.
	We refer the reader to \cite[Lemma 4.46]{DFL} for more details where a similar identity is established. 
	As a consequence of \eqref{h1-vanish}, we have
	\[
	  \int_{W}\left \vert d_A^h\zeta\right \vert^2 {\rm dvol}=2\int_{W}\left \vert (d_A^h\zeta)^+\right \vert^2{\rm dvol}.
	\]
	In particular, if $\zeta$ is in the kernel of $\mathcal D^h_A$, we have $d_A^h\zeta=0$. 
\end{proof}

\begin{lemma}\label{perturbed-ker-trivial}
	Suppose the inclusion of $N\smallsetminus U$ into $W\smallsetminus (\Bbb R\times U)$ induces a surjective map of $\pi_1$ and an isomorphism of $\pi_0$. 
	Let $A$ be an $h$-perturbed flat connection on $(W,c)$ such that the connections $\alpha\in \chi_h(N,w)$ and $\alpha'\in \chi_h(N',w')$ corresponding to the restriction of $A$ to the ends are 
	non-degenerate critical points of perturbed Chern--Simons functionals.
	Then $\ker(\mathcal D^h_A)$ is trivial for any $h$-perturbed flat connection $A$ on $(W,c)$.
\end{lemma}
\begin{proof}
	Lemma \ref{char-ker-coker} implies that any $\zeta$ in the kernel of $\mathcal D^h_A$ satisfies $d_A^h\zeta=0$. 
	By following the proof of \cite[Proposition 4.60]{DFL}, we can show that there is an $L^2_{k+1}$ section $\eta$ of $\su(E)$ such that $d_A\eta=\zeta$.
	To be a bit more specific, to define $\eta$ at a point $q\in W$, we take a path $\gamma:(-\infty,0]\to W$ such that $\gamma(0)=q$ and $\gamma(t)=(t,q_0)\in (-\infty,-1]\times N$
	for a fixed $q_0\in N$ and small enough values of $t$. In the case that $q=(s, q_0)\in \Bbb R\times U$, we require that $\gamma(t)\in \Bbb R\times \{q_0\}$, and in the case that $q$ is in the complement of $\Bbb R\times U$, we require that $\gamma(t)$ takes values in this complement for any $t$.
	Then $\eta(q)$ is given by integrating $\gamma^*\zeta$ over $(-\infty,0]$. (To make sense of this integration, 
	we use parallel transport with respect to the connection $\gamma^*A$ to trivialize the bundle.)

	To show that $\eta$ is well-defined, we need to show that this integral  
	depends only on $q$ and is independent of the choice of the path. This is straightforward in the case that $q\in \Bbb R\times U$, and we may assume that $q$ is 
	in the complement of $\Bbb R\times U$. This can be verified by showing that for any path $\gamma:(-\infty,\infty)\to W\smallsetminus (\Bbb R\times U)$, 
	satisfying $\gamma(t)=(t,q_0)\in (-\infty,-1]\times N$ for $t$ small enough and $\gamma(t)=(-t,q_0')\in (-\infty,-1]\times N$ for $t$ large enough, the integral 
	of $\gamma^*\zeta$ is trivial. For any $T\geq 1$ and using our assumption on the fundamental groups, we can find a compactly supported homotopy of 
	$\gamma$ in $W\smallsetminus (\Bbb R\times U)$ such that the image of $\gamma$ is contained in $(-\infty,-T]\times (N\smallsetminus U)$. 
	Then Stokes theorem and the triviality of $d_A\zeta$ in this region imply that any such homotopy does not change the value of the integral of $\gamma^*\zeta$. Now by letting 
	$T\to \infty$ and using the decay of $\zeta$ along the cylindrical ends, we can conclude that the integral of $\gamma^*\zeta$ is trivial. 
	
	Next, observe that
	\begin{align*}
		\int_{W} \left\vert \zeta\right \vert^2 {\rm dvol}&=\int_{W} \langle d_A\eta, \zeta \rangle\\
		&=\int_{W} \langle \eta, d_A^*\zeta \rangle=0,
	\end{align*}
	where the last identity follows from the assumption that $\zeta$ is in $\ker(\mathcal D^h_A)$, and hence $d_A^*\zeta$ vanishes. As a consequence of the above identity, $\zeta$ is trivial.
\end{proof}

The following lemma is the counterpart of Lemma \ref{char-ker-coker} for the cokernel of $\mathcal D^h_A$.
\begin{lemma}\label{char-coker}
	Suppose $A$ is an $h$-perturbed flat connection on $(W,c)$ that is asymptotic to non-degenerate perturbed projectively flat connections on the ends. 
	Then the cokernel of the ASD operator $\mathcal D^h_A$ is given by
	\begin{equation}\label{coker-lin-const}
		\coker(\mathcal D^h_A)=\{(0,\kappa)\in L^2_k(W, \su(E)\oplus \Lambda^+\otimes \su(E))\mid d_A^h\kappa=0\},
	\end{equation}
	where for any $\kappa\in L^2_{k}(W,\Lambda^2\otimes \su(E))$, the 3-form $d_A^h\kappa\in L^2_{k-1}(W,\Lambda^3\otimes \su(E))$ is defined as 
	\[
	  d_{A}^{h}\kappa :=d_{A}\kappa-*_3{\rm Hess}_{A_t}h((\iota_{\partial_t}\kappa)_t)\wedge dt.
	\]  
\end{lemma}

\begin{proof}
The non-degeneracy of $A$ at the ends is used to ensure that $\mathcal D^h_A$ is a Fredholm operator. When this holds, we can identify the cokernel of $\mathcal D^h_A$ with the kernel of its formal adjoint. The proof then follows by a similar argument as for Lemma \ref{char-ker-coker}. (See also the discussion preceding \cite[Lemma 4.57]{DFL}.)
\end{proof}

Let $(N,w)$, $K\subset N$, the 2-handle cobordism $(W,c)$ be as in Example \ref{example-cylinder-surgery-cob}, and $U$ be a codimension zero submanifold in the exterior of $K$ such that $K$ is null-homotopic in $N\smallsetminus U$.  There is an immersed sphere $\sigma:S^2\to W\smallsetminus (\Bbb R\times U)$ consisting of the union of the core  $\sigma_+$ of the 2-handle and a nodally immersed disk $\sigma_-$ in $(-\infty,-3]\times (N\smallsetminus U)$. Here we assume that the boundary of the core of the 2-handle $\sigma_+$ is contained in $\{-3\}\times N $. Similarly, the co-core of the 2-handle gives an embedded disk $\sigma_+'$  with boundary in $\{3\}\times N' $, and we glue to $\sigma_+'$ the subset $[3,\infty)\times \partial \sigma_+'$ of $[3,\infty)\times N'$ to obtain the embedded open disk $\sigma'$ in $W$. We write $p$ for the intersection point of $\sigma$ and $\sigma'$. Fix a cylinder function $h$ on $U$, and let $A$ be an $h$-perturbed flat connection on $(W,c)$. As an analogue of Lemma \ref{perturbed-ker-trivial}, we wish to obtain information about the cokernel of $\mathcal D^h_A$. The following discussion will be used in the next section to show that the cokernel of $\mathcal D^h_A$ under some additional constraints is trivial. 

Let $V$ be the pullback of the normal bundle of the above immersed sphere with respect to the map $\sigma$, and extend $\sigma$ to an immersion $\widetilde \sigma$ from $V$ to a tubular neighborhood of the immersed sphere in $W$, which we may assume is still contained in $W\smallsetminus (\Bbb R\times U)$. The pullback $\widetilde \sigma^* A$ is a projectively flat connection on $V$, and its adjoint defines a rank three local system $\mathcal L_A$ on $V$. Since this local system is necessarily trivial, we may fix an isomorphism between $H_c^2(V;\mathcal L_A)$ and the fiber $\mathcal L_p$ of $\mathcal L_A$ over $p$. This isomorphism maps $v \in \mathcal L_p$ into $\widetilde v\cdot \eta$, where $\widetilde v$ is the extension of $v$ into a parallel section of $\mathcal L_A$, and $\eta$ is a closed differential 2-form with compact support that integrates to $1$ over the fibers of $V$. We may push forward $\widetilde v\cdot \eta$ to $W$ to obtain 
\begin{equation}\label{phi-A-v}
	\phi_A(v)\in \Omega^2(W,\su(E))
\end{equation} 
that is supported in the chosen neighborhood of $\sigma$, and $d_A^h\phi_A(v)$ is trivial. Note that integrating $\phi_A(v)$ over the disk $\sigma'$ gives the element $v\in \mathcal L_p$.

\begin{lemma}\label{perturbed-coker-trivial}
	Let the perturbed projectively flat connection $A$ be asymptotic to non-degenerate perturbed projectively flat connections on the ends of $(W,c)$.
	Then for any $\kappa\in L^2_{k}(W,\Lambda^+\otimes \su(E))$ satisfying $d_{A}^{h}\kappa=0$, there are $v\in \mathcal L_p$ and 
	$\xi\in L^2_{k}(W,\Lambda^1\otimes \su(E))$ such 
	that $\kappa=\phi_A(v)+d_A^h \xi$. Furthermore, if $\kappa$ is nonzero, then we have
	\[
	  -\int_{W}\tr(\phi_A(v)\wedge \phi_A(v))>0.
	\]
	In particular, $v$ is non-zero.
\end{lemma}

 \begin{proof}
 	We follow a similar proof as in \cite[Proposition 4.63]{DFL}. First we make an observation about the case of cylinders. 
	Suppose $A_0$ is a perturbed flat connection on $(-\infty,t_0)\times N$ in the temporal gauge. In particular, $A_0$ is the pullback of a perturbed flat connection $B$. 
	Suppose $\kappa_0\in L^2_{k}((-\infty,t_0)\times N,\Lambda^2\otimes \su(E))$ satisfies
	$d_{A_0}^{h}\kappa_0=0$ and decays exponentially fast as $t\to -\infty$. Thus, if we 
	define
	\begin{equation}\label{xi0}
	  \xi_0(t,y)=\int_{-\infty}^{t}\iota_{\partial_t}\kappa_0(\tau, y),
	\end{equation}
	then $\xi_0$ decays exponentially fast and it belongs to $L^2_{k}((-\infty,t_0)\times N,\Lambda^1\otimes \su(E))$.
	Using the assumption $d_{A_0}^{h}\kappa_0=0$, one can see that $\kappa_0=d_{A_0}^h\xi_0$. Note that the assumption that $A_0$ is in 
	temporal gauge is not necessary because we may always apply a 
	suitable gauge transformation to $A_0$ (and $\kappa_0$) such that this assumption is satisfied. Similarly, if $A_1$ is a perturbed flat connection on $(t_0,\infty)\times N'$ 
	in temporal gauge,
	and $\kappa_1\in L^2_{k}((t_0,\infty)\times N',\Lambda^2\otimes \su(E))$ satisfies $d_{A_1}^{h}\kappa_1=0$ and decays exponentially fast as $t \to \infty$, then 
	\begin{equation}\label{xi1}
	  \xi_1(t,y')=-\int_{t}^{\infty}\iota_{\partial_t}\kappa_1(\tau, y')
	\end{equation}
	belongs to $L^2_{k}((t_0,\infty)\times N',\Lambda^1\otimes \su(E))$ and we have $\kappa_1=d_{A_1}^h\xi_1$.

	There is an improvement of the above observation which shall be useful for us. In the previous paragraph, if $\kappa_0$ vanishes on $(t_0-\epsilon,t_0)\times N$ for some 
	$\epsilon>0$, then we may assume that the restriction of $\xi_0$ on 
	$(t_0-\epsilon/2,t_0)\times N$ is trivial by the following modification. First we may again assume that $A_0$ is the pullback of a non-degenerate perturbed flat connection $B$
	on $(N,w)$. The restriction of \eqref{xi0} to $(t_0-\epsilon,t_0)\times N$ is the pullback of a
	1-form $z$ over $N$ in the kernel of 
	 $*d_{B}+{\rm Hess}_{B}h$. Our non-degeneracy assumption on $B$ implies that $z=d_B\eta$ for a section $\eta$ of $\su(E)$ over $N$. 
	 Let $\rho:\Bbb R \to \Bbb R$ be a bump function that equals $1$ on $t\geq t_0-3\epsilon/5$ and $0$ on $t\leq t_0-4\epsilon/5$.
	 Then $\xi_0':=\xi_0-d_{A_0}(\rho(t)\eta)$ is trivial on $(t_0-\epsilon/2,t_0)\times N$ and $\kappa_0=d_{A_0}^h\xi_0'$.	

	Now, we turn to the proof of the first claim of this lemma for $\kappa\in L^2_{k}(W,\Lambda^+\otimes \su(E))$ satisfying $d_{A}^{h}\kappa=0$. 
	The integral of $\kappa$ over $\sigma'$ can be regarded as an element of $\mathcal L_p$. Changing $\kappa$ by an element of the form $d_A^h \xi$ does not change this integral, and the integral of $\phi_A(v)$
	over $\sigma'$ is $v$. In particular, to prove the claim, it suffices to consider the case that the integral of $\kappa$ over $\sigma'$ is trivial. Next, we subtract elements of 
	the form $d_A^h \xi$ in three stages such that in each step the support of the resulting 2-form shrinks and finally it 
	becomes empty in the last stage. 
	
	For the first step, we apply the observation in the first paragraph of the proof to the case that $A_1$ and $\kappa_1$ are the restrictions of $A$ and $\kappa$ to 
	$(1,\infty)\times N'$ to obtain $\xi_1$. Since $A_1$ is asymptotic to a non-degenerate perturbed flat connection on this end, $\kappa_1$ is exponentially decaying on this end.
	Let $f:\Bbb R \to \Bbb R$ be a function that is equal to $1$ over $[2,\infty)$ and is supported in $(1,\infty)$. Extend $f\cdot \xi_1$ trivially to the rest of $W$.
	Then $\kappa-d_A^h(f\xi_1)$ vanishes on $[2,\infty)\times N'$. After replacing $\kappa$ with this 2-form, we may assume that the support of $\kappa$ is contained 
	in the complement of $(2,\infty)\times N'$. In particular, $\kappa$ restricted to $\sigma'$ is a compactly supported 2-form and the integral of this 2-form vanishes. 
	Since $A$ is flat in a small enough neighborhood $Z$ in $W\smallsetminus \((3,\infty)\times N'\)$ of the co-core 
	and $\kappa$ is $d_A$-closed in 
	this neighborhood, we may find an element  $\xi'\in L^2_{k}(W,\Lambda^1\otimes \su(E))$, which is compactly supported in $Z$ and $d_A^h\xi'$ agrees with $\kappa$
	in a neighborhood of the co-core. Thus, after replacing $\kappa$ with $\kappa-d_A^h\xi'$, we may assume that the support of $\kappa$ is contained in the interior of a subset 
	of $W$ that may be identified with $(-\infty, 3)\times N$ such that this identification maps $(-\infty, -3]\times N$ and $(-\infty, 3]\times U\subset W$ 
	to themselves by the identity map.
	After making this identification, we apply the observation in the second paragraph to the case that $A_0$ and $\kappa_0$ are the restrictions of $A$ and $\kappa$ to this region 
	to obtain $\xi_0$, which we extend trivially to the rest of $W$. Note that we may again argue that $\kappa_0$ is exponentially decaying on this end because $A_0$ 
	is asymptotic to a non-degenerate perturbed flat connection on $N$. In particular, we have $\kappa=d_{A}^h\xi_0$, which completes the proof of the first claim.

	Next, note that for any $\kappa \in L^2_k(W,  \Lambda^+\otimes \su(E))$, the expression $-\tr(\kappa\wedge \kappa)$ is non-negative, and the integral 
	of this expression over $W$ is positive unless $\kappa\equiv 0$. Thus, to prove the second claim, it suffices to show that for any 
	$\xi \in L^2_{k}(W,\Lambda^1\otimes \su(E))$ and $\mu\in L^2_{k}(W,\Lambda^2\otimes \su(E))$ satisfying $d_A^h \mu=0$, we have 
	\[
	  \int_{W}\tr(d_A^h \xi \wedge \mu )=0.
	\]
	This claim is similar to \eqref{h1-vanish}, and in the case that $h=0$, it is an immediate consequence of Stokes' theorem. 
	The more general case can be verified as in \cite[Lemma 4.46]{DFL}.
 \end{proof}

\subsection{Proofs of the main theorems}\label{main}

Now we are ready to complete the proof of Theorem \ref{thm:Boileau} for the remaining cases that the surgery coefficient is $\pm 1$. As we remarked in the introduction, it suffices to focus on the case that the surgery coefficient is $-1$: we assume that $K$ is a nullhomotopic knot in a 3-manifold $Y$ such that the Dehn surgery $Y_{-1}(K)$ is homeomorphic to $Y$ in an orientation preserving way. In particular, the core $K'$ of $Y_{-1}(K)$ is also nullhomotopic by Lemma \ref{lem:nullhomotopic-core}. Let $W:Y\to Y_{-1}(K)$ denote the elementary cobordism given by gluing a 2-handle along $K$. The following is a crucial step in the proof of Theorem \ref{thm:Boileau}.

\begin{theorem}\label{iso-cob-map}
	The cobordism map $I^\sharp(W):I^\sharp(Y)\to I^\sharp(Y_{-1}(K))$ is an isomorphism.
\end{theorem}

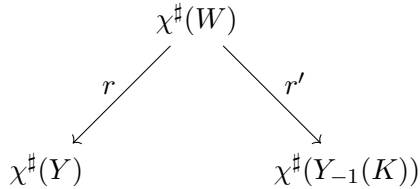
\begin{figure}[h]
\begin{center} 
  \begin{tikzpicture}
    \node (A) at (1,0) {$\chi^\sharp(W)$};
    \node (B) at (-1,-2) {$\chi^\sharp(Y)$};
    \node (C) at (3,-2) {$\chi^\sharp(Y_{-1}(K))$};
	\draw[->] (A) -- node[pos=0.6, above,xshift=-1pt] {$r$} (B);
	\draw[->] (A) -- node[pos=0.6, above, xshift=5pt] {$r'$} (C);
  \end{tikzpicture}
\end{center}
\caption{The correspondence associated to the character varieties}
  \label{fig:char-var-correspondence}
\end{figure}

We follow the notation at the end of subsection \ref{background-2} to denote $Y^\sharp$, $Y_{-1}^\sharp(K)$ and $W^\sharp$ for the admissible pairs associated to $Y $, $Y_{-1}(K)$ and the cobordism of admissible pairs determined by $W$. We also write $\chi^\sharp(W)$ for the character variety of $W^\sharp$, defined in the same way as in \eqref{char-ad-pair}, except that now we are considering a pair of a 4-manifold and a codimension two  submanifold. This character variety can be identified with energy zero solutions of the ASD equation (or equivalently projectively flat connections) on $W^\sharp$. As in \eqref{char-variety-sharp=rep-variety}, $\chi^\sharp(W)$ can be identified with two copies of the representation variety $R(W):={\rm Hom}(\pi_1(W),SU(2))$, and thus it is a real algebraic variety. The inclusion of $Y$ and $Y'$ in $W$ induces  restriction maps of character varieties as in Figure \ref{fig:char-var-correspondence}.

The assumption that $K$ (resp. $K'$) is null-homotopic implies that the inclusion of $Y$ (resp. $Y_{-1}(K)$) into $W$ induces an isomorphism of fundamental groups. In particular, the restriction maps $r$ and $r'$ determine isomorphisms of real algebraic schemes. Therefore, the map $r'r^{-1}:\chi^\sharp(Y) \to \chi^\sharp(Y_{-1}(K))$ is an isomorphism of real algebraic schemes, and it maps the Zariski tangent space $H^1(Y\#T^3,{\rm ad}_\alpha)$ at $\alpha\in \chi^\sharp(Y)$ to $H^1(Y_{-1}(K)\#T^3,{\rm ad}_{r'r^{-1}(\alpha)})$ by an isomorphism.

\begin{lemma}\label{lem:index}
	The index of the ASD operator $\mathcal D_A$ associated to a connection $A$ representing an element of $\chi^\sharp(W)$ and a metric on $W^\sharp$ with cylindrical ends 
	is zero. 
	In fact, the kernel and the cokernel of this operator are trivial. 
\end{lemma}
Note that we have not perturbed the ASD equation yet, and $\mathcal D_A$ is the unperturbed ASD operator. Because of this, we dropped $\pi$ from our notation for the ASD operator. Since the connection $A$ is not asymptotic to non-degenerate flat connections on $Y^\sharp$ and $Y_{-1}^\sharp(K)$ (except possibly when the dimension of the connected component of $\chi^\sharp(W)$  containing $[A]$ is zero), we need to modify the domain and the codomain of  $\mathcal D_A$ to obtain a well-behaved Fredholm operator. To do this, we use weighted Sobolev spaces where our function spaces are required to have exponential decay on the outgoing end and are allowed to have exponential growth on the incoming end, both for a small enough weight exponent $\delta$. 

We review this definition in some generality sufficient for our purposes. Suppose $(W,c):(N,w)\to (N',w')$ is a cobordism of admissible pairs. 
Let $L^2_{k,-\delta,\delta}(W)$ be the space of functions $f$ such that $\tau \cdot f \in L^2_{k}(W)$ where $\tau$ is a smooth function on $W$ which is equal to $e^{t\delta}$ on the half-cylinders $(-\infty,-1]\times N$ and $[1,\infty)\times N'$ with $t$ being the cylinder coordinates. More generally, for a vector bundle $V$ over $W$, we define the weighted Sobolev space $L^2_{k,-\delta,\delta}(W,V)$ as the space of sections $\sigma$ of $V$ such that $\tau\cdot \sigma \in L^2_{k}(W,V)$. For a connection $A$ on the bundle $E$ determined by $(W,c)$ asymptotic to connections $B$ on $(N,w)$ and $B'$ on $(N',w')$ over the ends, the ASD operator $(-d_A^*,d_A^+)$ determines a Fredholm operator 
\begin{equation}\label{ASD-weighted}
  \mathcal D_A:L^2_{k,-\delta,\delta}(W,\Lambda^1\otimes \su(E)) \to L^2_{k-1,-\delta,\delta}(W,(\Lambda^0\oplus \Lambda^+)\otimes \su(E)),
\end{equation}
whose index is independent of $\delta$, assuming it is a small enough positive real number, and Lemma \ref{lem:index} is an assertion about this index.

In the case that $A$ is a flat connection, the kernel of $\mathcal D_A$ has a characterization similar to Lemma \ref{char-ker-coker}. In fact a similar argument implies that 
\begin{equation}\label{ker-lin-const-weighted}
	\ker(\mathcal D_A)=\{\zeta \in L^2_{k,-\delta,\delta}(W,\Lambda^1\otimes \su(E))\mid d_A\zeta=0,\, d_A^*\zeta=0\}.
\end{equation}
This space is isomorphic to the relative cohomology group $H^1(W,N';{\rm ad}_{A})$ with twisted coefficients, where the local coefficient system ${\rm ad}_{A}$ is given by the flat connection $A$ on the vector bundle $\su(E)$. { Similarly to Lemma \ref{char-coker}, the cokernel of $\mathcal D_A$ may be identified with 
\begin{equation}\label{ker-lin-const-weighted}
	\coker(\mathcal D_A)=\{(0,\kappa)\in L^2_{k,\delta,-\delta}(W, \su(E)\oplus \Lambda^+\otimes \su(E))\mid d_A\kappa=0\}.
\end{equation}
This space is isomorphic to 
\begin{equation}\label{coker-DA}
	H^+(W;{\rm ad}_{A})\oplus \coker(H^1(W;{\rm ad}_{A}) \to H^1(N';{\rm ad}_{B'})), 
\end{equation} 
where $H^+(W;{\rm ad}_{A})$ is a maximal positive definite subspace of 
\[
  \im(H^2(W,N\sqcup N'; {\rm ad}_{A}) \to H^2(W; {\rm ad}_{A}))
\]
with respect to the intersection form. By definition, any element of \eqref{ker-lin-const-weighted} decays on the incoming end, and the first summand in \eqref{coker-DA} is given by the elements of \eqref{ker-lin-const-weighted} that also decay on the outgoing end. The second summand in \eqref{coker-DA} is determined by the possible limiting values of the elements of \eqref{ker-lin-const-weighted}. This can be regarded as a subspace of $H^1(N';{\rm ad}_{B'})$, which is the kernel of the Hessian of the Chern--Simons functional at $B'$.}

The above ASD index  is additive with respect to gluing connections; if $A$ is a connection on $(W,c):(N,w)\to (N',w')$ asymptotic to $B$ and $B'$ on the incoming and the outgoing ends and $A'$ is a connection on $(W',c'):(N',w')\to (N'',w'')$ asymptotic to $B'$ and $B''$ on the incoming and the outgoing ends, then we can glue $A$ and $A'$ to define a connection $A_\sharp$ on the composition $(W,c)\circ (W',c'):(N,w)\to (N'',w'')$ which is asymptotic to $B$ and $B''$ on the ends and we have
\begin{equation}\label{add-index-1}
  \ind(\mathcal D_{A_\sharp})=\ind(\mathcal D_{A})+\ind(\mathcal D_{A'}).
\end{equation}
Similarly, if $(W,c):(N,w)\to (N,w)$ is a cobordism from a pair to itself and $A$ is a connection on $(W,c)$ asymptotic to the same connection on the incoming and the outgoing ends, then we can glue the ends of $(W,c)$ to obtain a closed pair $(\overline W,\overline c)$ and $A$ induces a connection $\overline A$ on this pair. Then we have 
\begin{equation}\label{add-index-2}
  \ind(\mathcal D_{A})=\ind(\mathcal D_{\overline A}).
\end{equation}

Similarly to the non-degenerate case, there is a spectral flow characterization of the above ASD index in the case that $(W,c)$ is the cylinder corresponding to a pair $(N,w)$. Suppose $\{B(t)\}_{-1\leq t\leq 1}$ is a path of connections on $(N,w)$ from $B$ to $B'$. Then this path determines a connection $A$ on the product cobordism $\R\times (N,w)$ that is asymptotic to $B$ and $B'$ on the ends. The ASD index of the connection $A$ is equal to the spectral flow of the family of operators $L_{(B(t))}-\delta \cdot {\rm Id}$ for the same constant $\delta$ as above.

\begin{proof}[Proof of Lemma \ref{lem:index}]
	The above characterization of the kernel of the ASD operator for a flat connection implies that the kernel of $\mathcal D_A$ is isomorphic to 
	$H^1(W^\sharp,Y_{-1}^\sharp(K);{\rm ad}_{A})$, which is trivial. Thus, to complete the proof, it suffices to show that the index of $\mathcal D_A$ 
	is zero, which then implies the triviality of the cokernel of $\mathcal D_A$.
	
	Let $A$ be asymptotic to flat connections $B$ and $B'$ on the incoming and the outgoing ends of $W^\sharp$.
	Fix a diffeomorphism $\phi:Y\to Y_{-1}(K)$. Then $\phi^*B'$ is not necessarily isomorphic to $B$. However, we can find flat connections $A_0,\dots,A_k$ on $W^\sharp$ 
	for some $k$ such that the following holds:
	\begin{itemize}
		\item[(i)] $A_0=A$;
		\item[(ii)] if $A_i$ is asymptotic to $B_i$ and $B_i'$ on the incoming and the outgoing ends of $W^\sharp$, then $B_i=\phi^*B_{i-1}'$ 
		for any $1\leq i\leq k$.
		\item[(iii)] there is a path $\{B(t)\}_{0\leq t\leq 1}$ of flat connections on $Y^\sharp$ such that $B(0)=\phi^*B_{k}'$, $B(1)=B_0$, and the dimension of 
		$H^1(Y^\sharp, {\rm ad}_{B(t)})$ is independent of $t$.
	\end{itemize}
	
	{ To see why we may arrange for such a sequence of connections, first we consider the decomposition of $\chi^\sharp(Y)$ into the union of the following 
	subspaces:
	\[
	  \chi^\sharp(Y)^{(i)}=\{[B]\in \chi^\sharp(Y)\mid \dim(H^1(Y^\sharp, {\rm ad}_{B}))=i\}.
	\]
	This partition of $\chi^\sharp(Y)$ is the stratification of this space as a real algebraic variety into the subspaces where the dimension of Zariski tangent spaces is fixed.
	In particular, there are finitely many integers $i$ such that $\chi^\sharp(Y)^{(i)}$ is non-empty, and any such stratum has finitely many connected components; see \cite[Tag 05P8]{stacks-project}, \cite[Theorem 2.4.5]{bochnak2013real}. The pullback map $\phi^*: \chi^\sharp(Y_{-1}(K)) \to \chi^\sharp(Y)$ is an isomorphism of real algebraic varieties because 
	$\phi$ is a diffeomorphism. As we discussed above (just before the statement of Lemma~\ref{lem:index}), the maps $r'$, $r$ are also isomorphisms of real algebraic varieties. In particular, 
	$\phi^* \circ r' \circ r^{-1}$ is an automorphism of $\chi^\sharp(Y)$ as a real algebraic variety. Thus, it maps $\chi^\sharp(Y)^{(i)}$ to itself by a homeomorphism. 
	In particular, it induces a bijection of the set of connected components of $\chi^\sharp(Y)^{(i)}$. Thus there is $k$ such that 
	$\phi^*[B_k']=\(\phi^* \circ r' \circ r^{-1}\)^{k+1}([B_0])$ and $[B_0]$ belong to the same connected component of a stratum $\chi^\sharp(Y)^{(i)}$. This implies that there is 
	a piecewise smooth path from $\phi^*[B_k']$ to $[B_0]$ inside $\chi^\sharp(Y)^{(i)}$ because $\chi^\sharp(Y)^{(i)}$ is a semi-algebraic set; see \cite[Proposition 2.5.13]{bochnak2013real}. By reparametrizing this path, we may obtain a smooth path in $\mathcal B(E)$, which then can be lifted to a smooth path $\{B(t)\}_{0\leq t\leq 1}$ in 
	$\mathcal A(E)$ satisfying the requirements in (iii).

	The path $B(t)$ determines a connection $A_{k+1}$ on $\R\times Y^\sharp$ which has trivial topological energy and index. 
	The computation of the index follows from the characterization of the ASD index in terms of spectral flow and the fact that there is a small 
	enough positive constant $\delta$ such that any positive eigenvalue of the operator $L_{B(t)}$ is larger than $\delta$. {(Here we are using the assumption 
	that the dimension of the kernel of $L_{B(t)}$ is independent of $t$.)} 
	Furthermore, the observation in the first paragraph of the proof implies that the index of each $A_i$ with $0\leq i\leq k$ is non-positive. 
	
	The connections $A_0,\dots,A_{k+1}$ can be glued to each other in the cyclic order to define a connection $\bf A$ on the closed $4$-manifold $X$ 
	given by $k+1$ copies of $W^\sharp$,
	where we glue the incoming end of the $i^{\rm th}$ copy of $W^\sharp$ to the outgoing end of the $(i-1)^{\rm st}$ copy using the diffeomorphism $\phi$. 
	Since $W$ has Euler characteristic $1$ and signature $-1$, additivity of the signature and Euler characteristic with respect to gluing imply that the Euler 
	characteristic and the signature of $X$
	are respectively equal to $k+1$ and $-(k+1)$. The topological energy of $\bf A$, which is the sum of the topological energies of $A_i$, is equal to zero. Thus, the index formula for 
	the ASD operator on closed 4-manifolds implies that the index of $\bf A$ is zero. Moreover, the relations \eqref{add-index-1} and \eqref{add-index-2}
	assert that 
	\[
	  \ind(\mathcal D_{A_0})+\dots +\ind(\mathcal D_{A_{k+1}})=\ind(\mathcal D_{{\bf A}})=0.
	\]
	Since all the terms on the left hand side are non-positive, this identity implies that all the indices appearing in this formula are zero, completing the proof of the claim. }
\end{proof}

\begin{cor}
	Let $E$ be the $U(2)$-bundle associated to the admissible pair $W^\sharp$, and $v\neq 0$ be an element of the fiber of $\su(E)$ over a basepoint on the core of the 2-handle of 
	$W^\sharp$. Then we have
	\begin{equation}\label{neg-def-pairing}
	  -\int_{W^\sharp} \tr(\phi_A(v)\wedge \phi_A(v))<0.
	\end{equation}
	where $\phi_A(v)$ is constructed as in \eqref{phi-A-v}.
\end{cor}
\begin{proof}
	Lemma \ref{lem:index} and the identification of the cokernel of $\mathcal D_A$ with \eqref{coker-DA} imply that the space $H^+(W^\sharp;{\rm ad}_{A})$ is trivial. In particular, the intersection pairing on the image of the map 
\begin{equation}\label{relative-to-absolute}
  H^2(W^\sharp,Y^\sharp\sqcup Y_{-1}^\sharp(K);{\rm ad}_{A}) \to H^2(W^\sharp;{\rm ad}_{A})
\end{equation}
is negative definite. This map factors through $H^2(W^\sharp,Y^\sharp;{\rm ad}_{A})$, which is three-dimensional by excision. In fact, the forms $\phi_A(v)$, as $v$ varies, provide representatives for the image of the morphism in \eqref{relative-to-absolute}. To see this note that the co-core of $W^\sharp$ together with an immersed disk bounding $K'$ determines an immersed sphere in $W^\sharp$ such that the integral of $\phi_A(v)$ over this sphere is $v$. In particular, $\phi_A(v)$ is a non-trivial element of the image of the map in \eqref{relative-to-absolute} for a non-zero $v$. Since the intersection pairing on this image is negative definite, we obtain \eqref{neg-def-pairing}.
\end{proof}

\begin{remark}
	In the proof of Lemma \ref{lem:index}, it is essential to know that $Y_{-1}(K)$ is homeomorphic to $Y$ in an orientation preserving way so that 
	the glued up 4-manifold $X$ is oriented. This is the only step that we use this assumption, and for the remaining steps in the proof it suffices to only 
	assume that $Y_{-1}(K)$ and $Y$ have isomorphic fundamental groups. In Section~\ref{sec:example}, we discuss an example where 
	Lemma \ref{lem:index} does not hold when $Y_{-1}(K)$ and $Y$ are homeomorphic in an orientation reversing way. $\diamd$
\end{remark}

Next, we need to pick an appropriate perturbation of the Chern--Simons functional of the admissible pair $Y^\sharp$, which in turn induces a perturbation of the Chern--Simons functional for $Y_{-1}^\sharp(K)$ and a gradient-like perturbation of the ASD equation for $W^\sharp$. To find such a perturbation, we need the following topological lemma.

\begin{figure}
\begin{subfigure}{0.5\textwidth}
\labellist
\small\hair 2pt
\pinlabel {$K_0$} at 205 195
\pinlabel {$K_0$} at 405 195
\pinlabel {{\color{blue}$a_i$}} at 305 430
\pinlabel {$D_0$} at -20 400
\pinlabel {$D_0$} at 640 400
\endlabellist
\centering
\includegraphics[height=0.25\textheight, width=0.8\textwidth]{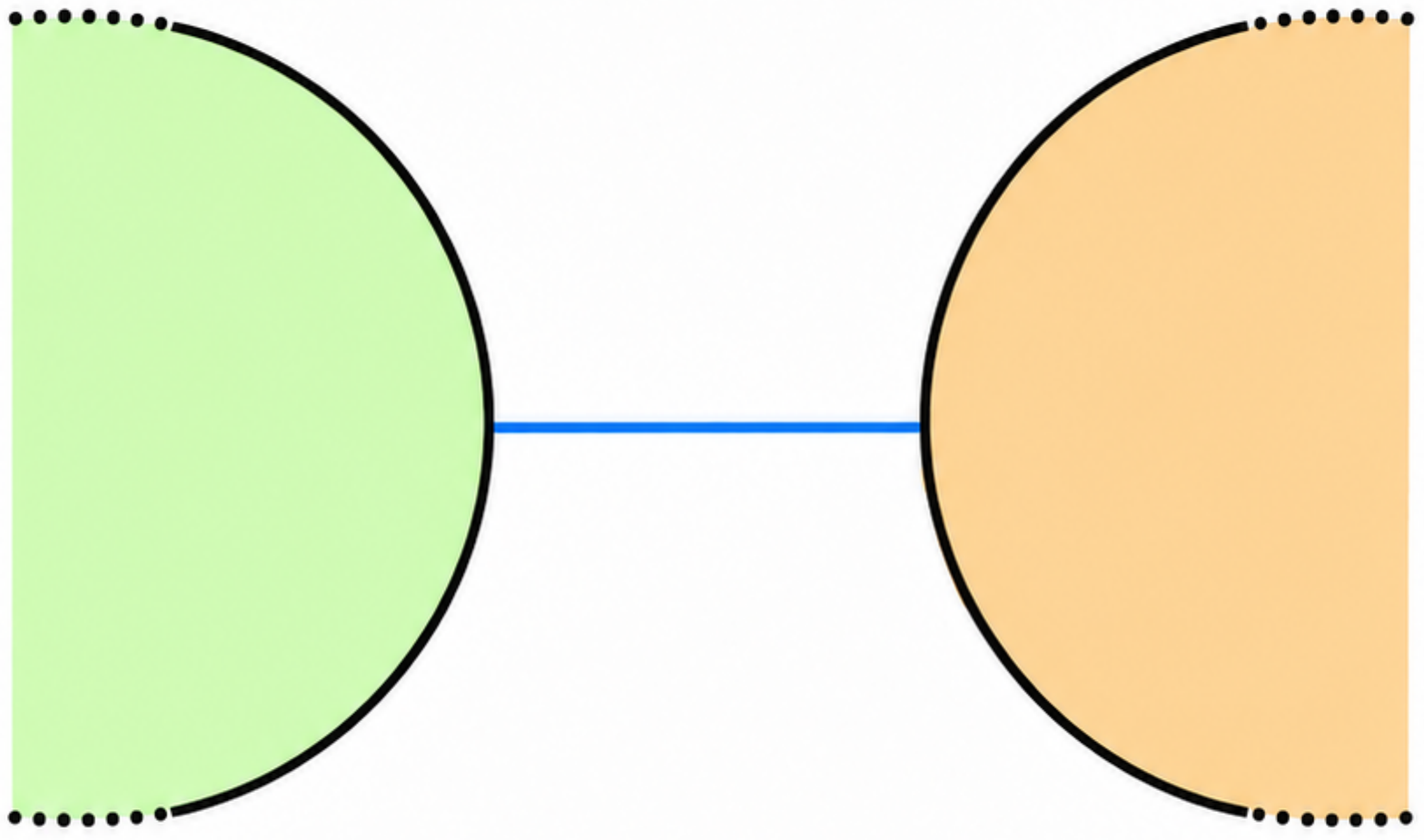}
\end{subfigure}
~
\begin{subfigure}{0.5\textwidth}
\labellist
\small\hair 2pt
\pinlabel {$K$} at 205 195
\pinlabel {$K$} at 405 195
\pinlabel {$D$} at -20 400
\pinlabel {$D$} at 640 400
\endlabellist
\centering
\includegraphics[height=0.25\textheight, width=0.8\textwidth]{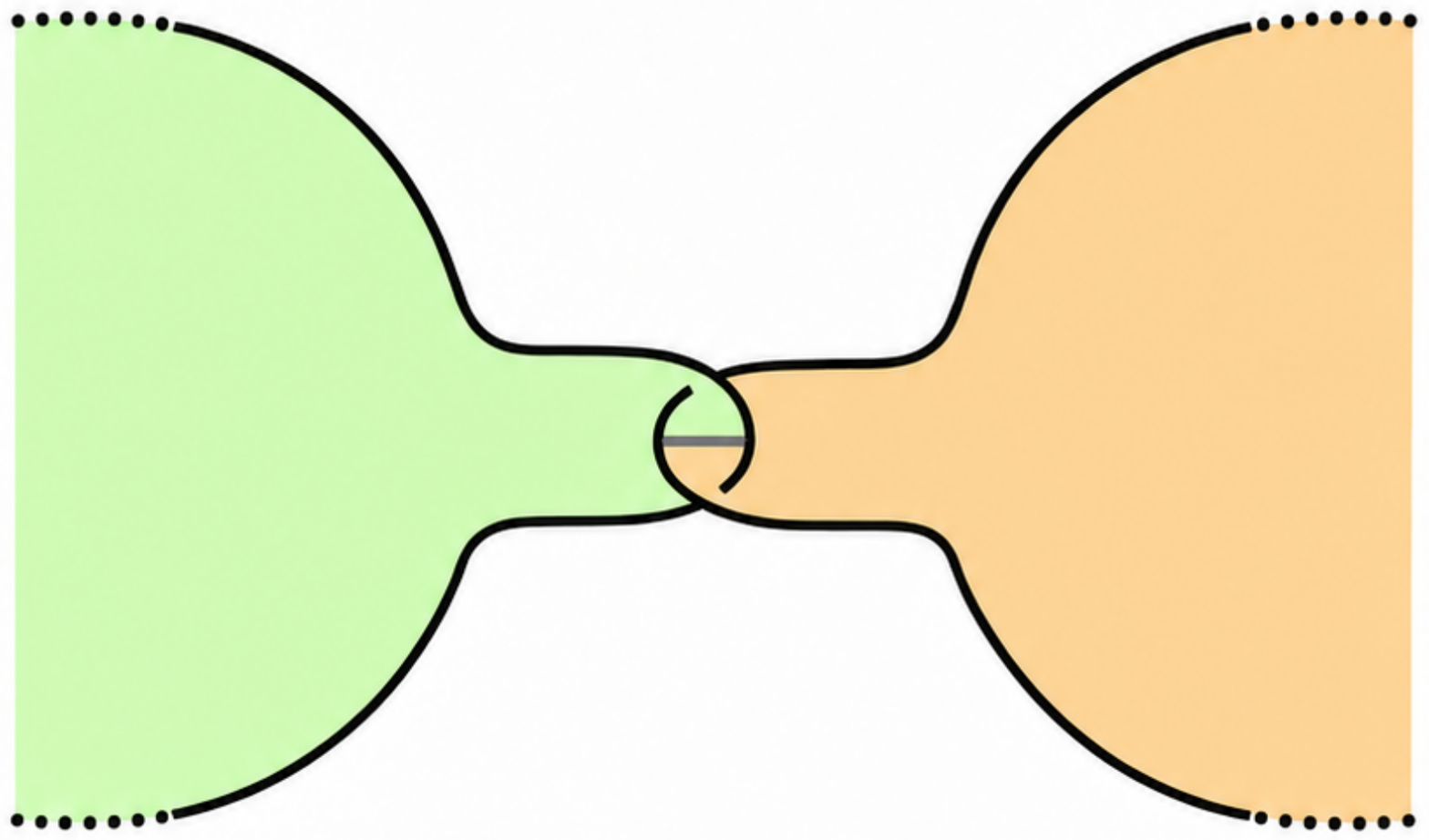}
\end{subfigure}
\caption{The left figure shows the unknot $K_0$ with one of the arcs and an embedded disk $D_0$ with boundary $K_0$, and the right figure shows how $D_0$ can be turned into an immersed disk with boundary $K$.}    \label{immersed-disk}
\end{figure}

\begin{lemma}\label{immersed-disc}
	Let $K$ be a nullhomotopic knot in a three-manifold $Y$.  Then, there exists a choice of an immersed disc $i:D \to Y$ such that the boundary of this disc is 
	$K$ and the inclusion of $Y - i(D)$ into $Y$ induces a surjection of fundamental groups.
\end{lemma}
\begin{proof}
	Since $K$ is nullhomotopic, it can be converted into the unknot by a sequence of $k$ crossing changes.  Alternatively, we can view $K$ as a sequence of $k$ 
	crossing changes on an unknot $K_0$ in $Y$.  Consider a sequence of disjoint arcs, $a_1,\ldots,a_k$ in $Y$ with boundary on $K_0$ describing the crossing changes. 
	Let also $D_0$ be a disc with boundary $K_0$. By putting the arcs $a_i$ in general position, we may assume that they intersect $D_0$ transversely.  
	We can find an immersed disc $D$ filling $K$ in any given neighborhood of $D_0 \cup a_1 \cup \ldots \cup a_k$ (see Figure \ref{immersed-disk}). We can fix the latter neighborhood such that 
	it is given by a 3-ball with several 1-handles attached (possibly more than $k$ if the $a_i$ intersect $D_0$ in the interior), i.e. a 3-dimensional handlebody.  
	Since $Y$ is obtained by attaching 2- and 3-handles to the complement of this handlebody, the result follows for this choice of $D$.
\end{proof}

{ The disk $D$ provided by Lemma \ref{immersed-disc} determines a disk in $Y \# T^3$, which is still denoted by $D$. By Lemma \ref{immersed-disc}, we may fix an open subspace $U$ of the complement of a neighborhood of this disk such that the inclusion of $U$ into $Y \# T^3$ induces a surjection of fundamental groups. This requirement on fundamental groups implies that for any distinct $\alpha_1$, $\alpha_2\in \chi^\sharp(Y)$, the restrictions of $\alpha_1$ and $\alpha_2$ to $U$ are not equal to each other. Similarly, for any $\alpha \in\chi^\sharp(Y)$ and any nonzero $\zeta\in H^1(Y\#T^3,{\rm ad}_{\alpha})$, the restriction of $\zeta$ to $U$ is not trivial in $H^1(U,{\rm ad}_{\alpha})$. Now using the arguments in \cite[Section 5.5.1]{donaldson-book}, we may find cylinder functions $h_1$, $h_2$, $\dots$, $h_n$ as in \eqref{special-cylinder} such that each $h_i$ is supported in $U$, the map \begin{equation}\label{family-pert}
	\underline h:=(h_1,\dots,h_n):\mathcal B(E) \to \Bbb R^n
\end{equation} 
restricted to $\chi^\sharp(Y)$ is injective, and the restriction of the derivative of $\underline h$ to $H^1(Y\#T^3,{\rm ad}_{\alpha})$, as a subspace of the tangent space of $\mathcal B(E)$ at $\alpha$, is injective for any $\alpha\in \chi^\sharp(Y)$.

Next, we use these cylinder functions to define a family version of $\chi^\sharp(Y)$. 

\begin{definition}
Let ${\bf X}^\sharp(Y)$ be the subspace of $\mathcal B(E)\times \Bbb R^n$ given by the tuples $([B],r_1,\dots,r_n)$ such that 
\begin{equation}\label{family-flat-eq}
  *F_B+(r_1\nabla_Bh_1+\dots+ r_n\nabla_Bh_n)=0.
\end{equation}
\end{definition}

Our assumptions on $h_i$ imply that any solution $B$ of this equation is a flat connection on the complement of $U$. Since the complement of $U$ contains a regular neighborhood of $K$ and $K$ is null-homotopic in this complement, any such connection is gauge equivalent to a connection that is trivial in the regular neighborhood of $K$.
There is an obvious projection map $\pi:{\bf X}^\sharp(Y) \to \Bbb R^n$ such that $\chi^\sharp(Y)=\pi^{-1}(0)$, and this projection map is proper by Uhlenbeck compactness theorem. For any positive real number $\epsilon$, define 
\[
  {\bf X}_\epsilon^\sharp(Y)=\{([B],r_1,\dots,r_n)\in {\bf X}^\sharp(Y)\mid |r_i|<\epsilon \text{ for all $i$}\}.
\] 
Moreover, for any cylinder function $h=r_1h_1+\dots+r_nh_n$, we define $\chi_h^\sharp(Y)$ as $\pi^{-1}(r_1,\dots,r_n)$.

The space ${\bf X}^\sharp(Y)$ is cut out transversely at any point of the subspace $\chi^\sharp(Y)$. To see this, note that the linearization of \eqref{family-flat-eq} at $([B],0)$ is given by 
\begin{align}\label{linearized}
  (\beta,\underline s)\in L^2_k(&Y \# T^3,\Lambda^1\otimes \su(E))\times \Bbb R^n \xrightarrow{}*d_B\beta+M(\underline s) \in {\rm ker}(d_B^*).
\end{align}
Here $M(\underline s):=s_1\nabla_Bh_1+\dots+ s_n\nabla_Bh_n$ where $\underline s=(s_1,\dots,s_n)$. As in Section \ref{background-2}, the cokernel of this operator may be identified with the cokernel of the operator 
\[
  \mathfrak L_B: L^2_k(Y \# T^3,(\Lambda^0\oplus \Lambda^1)\otimes \su(E))\times \Bbb R^n \to L^2_{k-1}(Y \# T^3,(\Lambda^0\oplus \Lambda^1)\otimes \su(E))
\]
that maps $(\xi,\beta,\underline s)$ to $L_B(\xi,\beta)+(0,M(\underline s))$. The cokernel of $\mathfrak L_B$ is given by the elements of the cokernel of $L_B$ that are also orthogonal to $\nabla_Bh_i$ for any $i$. Since $L_B$ is self-adjoint, any such element belongs to $H^1(Y\#T^3,{\rm ad}_{[B]})$ and is in the kernel of the derivative of \eqref{family-pert} at $[B]$. Our assumption on $\underline h$ implies that any such element is trivial. Consequently, if $\epsilon$ is small enough, \eqref{family-flat-eq} is cut out transversely at the elements of ${\bf X}_\epsilon^\sharp(Y)$. Therefore, it is a smooth manifold of dimension $\ind(L_B)+n=n$. 

The cylinder functions $h_i$ are supported in the complement of a tubular neighborhood of $K$. In particular, they induce perturbations of the Chern--Simons function for $Y_{-1}^\sharp(K)$, which we denote with the same notation. In particular, we obtain the analogue of \eqref{family-pert} as
\begin{equation}\label{family-pert-Y'}
	\underline h':=(h_1,\dots,h_n):\mathcal B(E') \to \Bbb R^n
\end{equation} 
where $E'$ denotes the $U(2)$-bundle associated to the admissible pair $Y_{-1}^\sharp(K)$. We also form 
\[
  {\bf X}^\sharp(Y_{-1}(K)), \hspace{1cm}{\bf X}_\epsilon^\sharp(Y_{-1}(K)), \hspace{1cm}\chi_h^\sharp(Y_{-1}(K))
\] 
as in the case of $Y$. We  write $\pi'$ for the counterpart of the map $\pi$. For any $\alpha\in \chi^\sharp(Y)$, we have $\underline h(\alpha)=\underline h'r'r^{-1}(\alpha)$ because the restrictions of $\alpha$ and $r'r^{-1}(\alpha)$ to $U$ agree. Similarly, the derivatives of $\underline h$ and $\underline h'r'r^{-1}$ restricted to $H^1(Y\#T^3,{\rm ad}_{\alpha})$ agree. In particular, the restriction of $\underline h'$ to $\chi^\sharp(Y_{-1}(K))$ is injective, and the restriction of the derivative of $\underline h'$ to any $H^1(Y_{-1}(K)\# T^3, {\rm ad}_{\alpha'})$ is injective for any $\alpha'\in \chi^\sharp(Y_{-1}(K))$. As in the above paragraph, we conclude that if $\epsilon$ is small enough, then ${\bf X}_\epsilon^\sharp(Y_{-1}(K))$ is cut out transversely, and it is a smooth manifold of dimension $n$.

Let $[B]\in \mathcal B(E)$ be a connection such that its restriction to the fixed tubular neighborhood of $K$ is gauge equivalent to the trivial connection. Then the restriction of $[B]$ to the exterior of $K$ is trivial on the boundary. This restriction may be also regarded as the gauge equivalence class of a connection on the exterior of $K'$ in $Y^\sharp_{-1}(K)$, that can be extended by the gauge equivalence class of the trivial connection to the tubular neighborhood of $K'$. This gives a well-defined gauge equivalence class of a connection on $Y^\sharp_{-1}(K)$ that depends only on $[B]$ and is denoted by $\psi([B])$. In particular, for any $([B],t_1,\dots,t_n)\in {\bf X}^\sharp(Y)$, $\psi([B])$ is defined and $(\psi [B],t_1,\dots,t_n)\in {\bf X}^\sharp(Y_{-1}(K))$ because the supports of the cylinder functions $h_i$ are in the complement of the fixed tubular neighborhood of $K$. It is straightforward to see that this defines a smooth map 
\[
  \Psi:{\bf X}_\epsilon^\sharp(Y) \to {\bf X}_\epsilon^\sharp(Y_{-1}(K)), \hspace{2cm} \Psi([B],t_1,\dots,t_n)=(\psi [B],t_1,\dots,t_n)
\]
that satisfies $\pi'\Psi=\pi$, the restriction of $\Psi$ to $\chi^\sharp(Y)$ is equal to $r'r^{-1}$, and the restriction of the derivative of $\Psi$ to $H^1(Y\#T^3,{\rm ad}_{\alpha})$ with $\alpha\in \chi^\sharp(Y)$ agrees with the homomorphism induced by $r'r^{-1}$ on this Zariski tangent space. In particular, the derivative of $\Psi$ restricted to such tangent spaces is injective.

\begin{lemma}\label{chi-h-chi'-h-smae-size}
	There is an arbitrarily small $h$ such that $\chi_h^\sharp(Y)$ and $\chi_h^\sharp(Y_{-1}(K))$ are non-degenerate and  have the same number of elements. 
\end{lemma}
\begin{proof}
	First we show that the map $\Psi$ is a local diffeomorphism at the points of $\chi^\sharp(Y)$. It suffices to show that at any point $([B],0)\in {\bf X}_\epsilon^\sharp(Y)$, the derivative of 
	$\Psi$ is injective. If $(\beta,\underline s)\in T_{([B],0)}{\bf X}_\epsilon^\sharp(Y)$ is in the kernel of the derivative $\Psi_*$ of $\Psi$, then we have 
	\begin{align*}
	  \underline s&=\pi_*(\beta,\underline s)\\
	  &=\pi'_*\Psi_*(\beta,\underline s)=0.
	\end{align*}
	This implies that this element of the kernel of $\Psi_*$ belongs to $H^1(Y\#T^3,{\rm ad}_{B})$. As we discussed above, the restriction of $\Psi_*$ to this Zariski 
	tangent space is injective, which gives the claim. Since $\pi$ and $\pi'$ are proper, $\pi'\Psi=\pi$, the map $\Psi$ is a local diffeomorphism at the points of $\pi^{-1}(0)$, and $\Psi$ maps 
	$\pi^{-1}(0)$ to $\pi'^{-1}(0)$ by a bijection, we have that $\Psi$ induces a diffeomorphism from ${\bf X}_\epsilon^\sharp(Y) $ to ${\bf X}_\epsilon^\sharp(Y_{-1}(K))$ after 
	possibly making $\epsilon$ smaller. In particular, for $h$ small enough, $\Psi$ maps $\chi_h^\sharp(Y)=\pi^{-1}(h)$ to $\chi_h^\sharp(Y_{-1}(K))=\pi'^{-1}(h)$ by a bijection.
	By Sard's theorem, we may pick such $h$ that is a regular value of $\pi$, $\pi'$ and it is arbitrarily small. Then, $\chi_h^\sharp(Y)$ and $\chi_h^\sharp(Y_{-1}(K))$
	are non-degenerate.
\end{proof}

Now, let $h$ be chosen such that it satisfies the claim in Lemma \ref{chi-h-chi'-h-smae-size}. In particular, $\chi_h^\sharp(Y)$ and $\chi_h^\sharp(Y_{-1}(K))$ consist of a finite set of non-degenerate critical points and have the same size. Fix a Riemannian metric on $W^\sharp$ with cylindrical ends as in Example \ref{example-cylinder-surgery-cob} such that it contains an isometric embedding of $\R\times U$ with product metric for the open set $U$ constructed above. In particular, the cylinder function $h=r_1h_1+\dots+r_nh_n$ induces a gradient like perturbation of the ASD equation on $W^\sharp$. We write $\chi_h^\sharp(W)$ for the space of all $h$-perturbed flat connections on $W^\sharp$. As in Figure \ref{fig:char-var-correspondence}, restriction maps to the boundary determine the following correspondence: 
\begin{center} 
  \begin{tikzpicture}
    \node (A) at (1,0) {$\chi_h^\sharp(W)$};
    \node (B) at (-1,-2) {$\chi_h^\sharp(Y)$};
    \node (C) at (3,-2) {$\chi_h^\sharp(Y_{-1}(K))$};
	\draw[->] (A) -- node[pos=0.6, above,xshift=-1pt] {$r_h$} (B);
	\draw[->] (A) -- node[pos=0.6, above, xshift=5pt] {$r_h'$} (C);
  \end{tikzpicture}
\end{center}
As it is explained in Example \ref{elem-cob-almost-flat}, we can identify $\chi_h^\sharp(W)$ as the space of all gauge equivalence classes of connections in $\chi_h^\sharp(Y)$ whose holonomy along $K$ is trivial and the restriction map $r_h$ is given by inclusion. Similarly, $\chi_h^\sharp(W)$ is given by the set of elements of $\chi_h^\sharp(Y_{-1}(K))$ whose holonomy along $K'$ is trivial and  $r_h'$ is the  inclusion map. In particular, $r_h$ is a bijection and $r_h'$ is injective. Since $\chi_h^\sharp(W)$ and $\chi_h^\sharp(Y_{-1}(K))$ have the same size (by Lemma~\ref{chi-h-chi'-h-smae-size}), we conclude that  $r_h'$ is also a bijection.

}

\begin{prop}\label{regular-almost-flat}
	If $h$ is small enough, then for any $h$-perturbed flat connection $A$, the ASD operator $\mathcal D^h_A$ has trivial kernel and cokernel. In particular, it has index zero.
\end{prop}

\begin{proof}
	The claim about the kernel follows from Lemma \ref{perturbed-ker-trivial} without making any assumption on the size of $h$. 	
	Next, we show that if $h$ is small enough, then the cokernel of $\mathcal D^h_A$ is trivial. If this claim does not hold, then there is a sequence $\{A_i\}$ with 
	$[A_i]\in \chi_{h_i}^\sharp(W)$ such that $h_i \to 0$ and the cokernel of $\mathcal D^{h_i}_{A_i}$ has a non-trivial element. By Lemma \ref{perturbed-coker-trivial}, 
	for any $i$, there is a non-zero vector $v_i$ in the fiber over a fixed basepoint on the 2-handle in $W^\sharp$ such that 
	\begin{equation}\label{ineq-vi}
	  -\int_{W}\tr(\phi_{A_i}(v_i)\wedge \phi_{A_i}(v_i))>0.
	\end{equation}
	Without loss of generality, we may assume that $v_i$ has norm $1$. By Uhlenbeck compactness, the sequence of connections $A_i$, after passing to a subsequence 
	and applying gauge transformations, is $C^\infty$-convergent on compact subspaces of $W^\sharp$ to a connection $A_\infty$ representing an element of $\chi^\sharp(W)$.
	We may also assume that $v_i$ is convergent to a unit vector $v_\infty$. Since $\phi_{A_i}(v_i)$ depends on $A_i$ on a fixed compact subspace of 
	$W^\sharp$, we conclude that $\phi_{A_i}(v_i) \to \phi_{A_\infty}(v_\infty)$. Thus \eqref{ineq-vi} implies that 
	\[
	   -\int_{W}\tr(\phi_{A_\infty}(v_\infty)\wedge \phi_{A_\infty}(v_\infty))\geq 0,
	\]
	which contradicts  \eqref{neg-def-pairing}. 
\end{proof}

\begin{proof}[Proof of Theorem \ref{iso-cob-map}]
 	 First we prove this theorem under some simplifying assumptions. We have already arranged the cylinder 
	function $h$ such that the perturbation of the Chern--Simons functionals of $Y^\sharp$ and  $Y_{-1}^\sharp(K)$ by $h$ have non-degenerate
	critical points. We temporarily assume that these perturbations can be used to form the instanton Floer complexes for $Y^\sharp$ and  $Y_{-1}^\sharp(K)$. That is to say, 
	the moduli spaces $\breve M_h(Y^\sharp;\alpha,\beta)_d$ and $\breve M_h(Y_{-1}^\sharp(K);\alpha',\beta')_d$ for $d\leq 1$ are regular. We also assume that the 
	moduli spaces $M_h(W^\sharp,\alpha,\alpha')_d$ of the perturbed ASD equation with $d\leq 1$ are regular where the perturbation is given by the gradient-like perturbation of 
	the ASD equation associated to $h$. Proposition \ref{regular-almost-flat} implies that this assumption holds for $h$-perturbed flat connections over $W^\sharp$. However, it 
	does not necessarily hold for other elements of these moduli spaces. 
	
	To show that the homomorphism $I(W^\sharp)$ associated to $W^\sharp$ is an isomorphism, we show that the corresponding chain map 
	$C(W^\sharp):C_*(Y^\sharp)\to C_*(Y_{-1}^\sharp(K))$ is an isomorphism. The Floer complexes $C_*(Y^\sharp)$ and $C_*(Y_{-1}^\sharp(K))$ are vector spaces generated
	by the elements of $\chi_h^\sharp(Y)$ and $\chi_h^\sharp(Y_{-1}(K))$. In the following, we write $J$ for the bijection $r_h'\circ r_h^{-1}$ and the induced isomorphism 
	$C_*(Y^\sharp)\to C_*(Y_{-1}^\sharp(K))$. 
	Lift the relative grading on $C_*(Y^\sharp)$ into an absolute grading in an arbitrary way and 
	then lift the grading on $C_*(Y_{-1}^\sharp(K))$ in a way that $C(W^\sharp)$ preserves the grading. Then Proposition  \ref{regular-almost-flat} implies that 
	for any generator $\alpha\in \chi_h^\sharp(Y)$ of $C_*(Y^\sharp)$, the corresponding generator $J(\alpha)$ of $\chi_h^\sharp(Y_{-1}(K))$ has the same 
	grading as $\alpha$. 
		
	We define an ordering on the elements of $\chi_h^\sharp(Y)$ with the same grading.  
	Let $\alpha$ and $\beta$ be elements of $\chi_h^\sharp(Y)$ whose gradings are equal to each other. Then there is a connection 
	$A_0$ on $\bbR\times Y^\sharp$ with ASD index zero whose restrictions to $(-\infty,-1]\times Y^\sharp$ and $[1,\infty)\times Y^\sharp$ 
	are respectively pullbacks of representatives 
	for $\alpha$ and $\beta$. For any other such connection $A_0'$, there is an element $u$ of the gauge
	group associated to the admissible pair $Y^\sharp$ such that $A_0'-u^*A_0$ is compactly supported. 
	In particular, \eqref{top-energy-h-rewrite} implies that $\mathcal E_h(A_0)$, defined in \eqref{top-energy-h} 
	with $W$ being the product cobordism associated to $Y^\sharp$, depends only on $\alpha$ and $\beta$, and we 
	write $\alpha \succeq \beta$ (resp. $\alpha \succ \beta$) if $\mathcal E_h(A_0)\leq 0$ (resp. $\mathcal E_h(A_0)< 0$).

	For connections $\alpha$, $\beta$ with 
	$\alpha \succeq \beta$, 
	let $\beta'=J(\beta)$ and $A$ be a connection over $W^\sharp$ 
	representing an element of the moduli space $M_h(W^\sharp,\alpha,\beta')_0$. Then Lemma \ref{pos-top-energy-sol} implies that
	$\mathcal E_h(A)\geq 0$. Now, let $A_0$ be as above and $A'$ be the connection obtained by gluing $A_0$ over $\R\times Y^\sharp$
	and the $h$-perturbed flat connection $r_h^{-1}(\beta)$. Proposition \ref{regular-almost-flat} and additivity of indices with respect to gluing imply 
	that the ASD operator associated to $A'$ has index zero. Since the connections $A$ and $A'$ on $W^\sharp$ have the same indices and 
	are asymptotic to the same connections on the ends, they have the same topological energy. Thus, \eqref{top-energy-h-rewrite} implies that 
	$\mathcal E_h(A)$ and $\mathcal E_h(A')$ are equal to each other. We can use this to see
	\begin{align*}
	  \mathcal E_h(A_0)&=\mathcal E_h(A_0)+\mathcal E_h(r_h^{-1}(\beta))\\
	  &=\mathcal E_h(A')=\mathcal E_h(A)\geq 0.
	\end{align*}
	This is a contradiction unless $\mathcal E_h(A)=0$, where another application of Lemma \ref{pos-top-energy-sol} shows that $\beta'=J(\alpha)$ or equivalently $\alpha=\beta$.
	In summary the moduli spaces $M_h(W^\sharp,\alpha,\beta')_0$, which are used to define the chain map $C(W^\sharp)$, are non-empty only if $\beta \succ \alpha$ or 
	$\alpha=\beta$. In the latter case, the moduli space $M_h(W^\sharp,\alpha,\beta')_0$ contains exactly one element which is $h$-perturbed flat. 
	Thus, for any $\alpha \in \chi_h^\sharp(Y)$, the element $C(W^\sharp)(\alpha)$ of $C_*(Y_{-1}^\sharp(K))$ is equal to $\pm J(\alpha)$ plus a linear combination of terms with higher order with respect to 
	$\succeq$. In particular, $C(W^\sharp)$ is an isomorphism.
	
	Next, we turn to the general case, where we might need additional perturbation terms for the Chern--Simons functionals and the ASD equation. 
	First we perturb the cylinder function $h$ further to guarantee that the instanton homology groups of $Y^\sharp $ and $Y_{-1}^\sharp(K)$ are defined. 
	Using the arguments in  \cite[Section 5.5.1]{donaldson-book}, we can find arbitrarily 
	small cylinder functions $g$ and $g'$ associated to $Y^\sharp $ and $Y_{-1}^\sharp(K)$ such that the following properties hold:
	\begin{itemize}
		\item[(i)] $g$ (resp. $g'$) is trivial in a neighborhood of $\chi_h^\sharp(Y)$ (resp. $\chi_h^\sharp(Y_{-1}(K))$), and the critical points of 
		the perturbation of the Chern--Simons functional of $Y^\sharp$ (resp. $Y_{-1}^\sharp(K)$) by $h+g$ (resp. $h+g'$) are $\chi_h^\sharp(Y)$ 
		(resp. $\chi_h^\sharp(Y_{-1}(K))$); 
		\item[(ii)] the moduli spaces $\breve M_{h+g}(Y^\sharp;\alpha,\beta)_d$ and $\breve M_{h+g'}(Y_{-1}^\sharp(K);\alpha',\beta')_d$ for $d\leq 1$ are regular. 
	\end{itemize}	
	
	There is also a secondary perturbation term $\omega$ for the ASD equation on $W^\sharp$ such that the moduli spaces $M_\pi(W^\sharp,\alpha,\alpha')_d$ of solutions to the 	equation 
	\begin{equation}\label{perturbed-ASD-W-sharp}
	  F_A^++(*_3\nabla_{A_t}h)^++\rho_-(t)(*_3\nabla_{A_t}g)^++\rho_+(t)(*_3\nabla_{A_t'}g')^++\omega(A)=0
	\end{equation}
	with $d\leq 1$, $\alpha\in \chi_h^\sharp(Y)$ and $\alpha'\in \chi_h^\sharp(Y_{-1}(K))$ consist of regular elements. 
	Here $\pi=(\omega,h+g,h+g')$, and we can assume that $\omega$ is as small as we wish.
	Using Uhlenbeck--Floer compactness and the observation in the simplified case, we see that if $g$, $g'$ and $\omega$ are small enough and 
	$\beta'=J(\beta)$, then the moduli spaces 
	$M_\pi(W^\sharp,\alpha,\beta')_0$ with $\alpha \succeq \beta$ are empty unless $\alpha=\beta$.
	
	Next, we show that $M_\pi(W^\sharp,\alpha,\alpha')_0$ with $\alpha'=J(\alpha)$ consists of exactly one element if $g$, $g'$ and $\omega$ are small enough.
	 Let $A^\alpha=r_h^{-1}(\alpha)$. Proposition \ref{regular-almost-flat} implies that
	\[
	  \mathcal D_{A^\alpha}^h:L^2_{k}(W,\Lambda^1\otimes \su(E)) \to L^2_{k-1}(W,(\Lambda^0\oplus \Lambda^+)\otimes \su(E))
	\]
	is an invertible operator. In particular, there is an operator 
	\[
	  Q:L^2_{k-1}(W,\Lambda^+\otimes \su(E)) \to \ker(d_{A^\alpha}^*)\subseteq L^2_{k}(W,\Lambda^1\otimes \su(E))
	\] 
	such that 
	\[
	  Q\circ d_{A^\alpha}^{+,h}(a)=a,\hspace{1cm}d_{A^\alpha}^{+,h}\circ Q(z)=z,\] 
	where $d_{A^\alpha}^{+,h}:=d_{A^\alpha}^++(*_3{\rm Hess}_{A^\alpha_t}h)^+$, $a\in\ker(d_{A^\alpha}^*)$ and $z\in L^2_{k-1}(W,\Lambda^+\otimes \su(E))$. 
	Using the Coulomb gauge slice, any connection in a small neighborhood of the gauge equivalence class of $A^\alpha$ has the form $A^\alpha+a$ for 
	$a\in \ker(d_{A^\alpha}^*)$. For any connection of this form, equation \eqref{perturbed-ASD-W-sharp} has the form
	\begin{equation}\label{perturbed-ASD-W-sharp-2}
	  d_{A^\alpha}^{+,h}(a)=P(a),
	\end{equation}
	where 
	\[P(a)=-(a\wedge a)^+-q_h(a)-\rho_-(t)(*_3\nabla_{\alpha+a_t}g)^+-\rho_+(t)(*_3\nabla_{\alpha'+a_t'}g')^+-\omega(A^\alpha+a),\] 
	with $$q_h(a)=[*_3(\nabla_{A^\alpha_t+a_t}h-\nabla_{A^\alpha_t}h-{\rm Hess}_{A^\alpha_t}h(a_t))]^+.$$
	By setting $a=Q(z)$ for $z\in L^2_{k-1}(W,\Lambda^+\otimes \su(E))$, we see that \eqref{perturbed-ASD-W-sharp-2} can be written as the following fixed point equation:
	\begin{equation}\label{fixed-point}
	  PQ(z)=z.
	\end{equation}
	It is straightforward to see that if $g$, $g'$ and $\omega$ are small enough, then $PQ(z)$ defines a contraction on a neighborhood of the origin in 
	$L^2_{k-1}(W,\Lambda^+\otimes \su(E))$, and hence  \eqref{fixed-point} has a unique fixed point. A compactness argument as in the previous case shows that there is 
	no solution outside of this neighborhood. This means that $M_\pi(W^\sharp,\alpha,\alpha')_0$ with $\alpha'=J(\alpha)$ has exactly one point. Therefore, as in the simplified case,
	$C(W^\sharp)$ is an isomorphism.
\end{proof}

\begin{proof}[Proof of Theorem~\ref{thm:Boileau}]
	Let $K\subset Y$ be a non-trivial null-homotopic knot such that $Y_{m/n}(K)$ is homeomorphic to $Y$ in an orientation preserving way. By Theorem \ref{thm:|n|>1}, 
	we may assume that $m/n=-1$.
	The surgery exact triangle of \cite{Scaduto} applied to the knot $K$ in the admissible pair $Y^\sharp=(Y\#T^3,S^1)$ implies that there is an exact triangle of instanton Floer homology groups 
	of the form
	\begin{equation}\label{exact-triangle}
		\cdots \xrightarrow{ } I^\sharp(Y)\xrightarrow{I^\sharp(W)} I^\sharp(Y_{-1}(K))\xrightarrow{ }I(Y_0(K)\#T^3, w_0 \cup S^1) \xrightarrow{ } I^\sharp(Y)\xrightarrow{I^\sharp(W)} 
		\cdots
	\end{equation}
	where $w_0$ is the embedded 1-manifold given by the core of surgery in $Y_0(K)$. Since $I^\sharp(W)$ is an isomorphism by Theorem \ref{iso-cob-map},
	\eqref{exact-triangle} implies that $I(Y_0(K)\#T^3, w_0 \cup S^1)$ is trivial. Thus, we can apply Theorem \ref{non-triviality-I-detailed} to conclude that $K$ is the trivial knot,
	which is a contradiction.
\end{proof}

\begin{proof}[Proof of Theorem~\ref{thm:Boileau-2}] Suppose $f$ is an (orientation-preserving) homeomorphism from the complement of $K'$ to the complement of $K$. Extending $f$ over the filling gives a homeomorphism from $Y$ to some surgery on $K \subset Y$. Lemma~\ref{lem:nullhomotopic-core} shows that $K'$ is nullhomotopic.

If $K$ and $K'$ are both trivial, the claim follows immediately. If one of them is non-trivial, then without loss of generality we can assume it is $K$. By applying Theorem~\ref{thm:Boileau} we deduce that the surgery is trivial, and the conclusion follows. 
\end{proof}

\begin{proof}[Proof of Theorem~\ref{thm:Boileau-3}]  Suppose $K \subset Y$ and $K' \subset Y'$ are null-homotopic and have homeomorphic exteriors. By changing the orientation of $Y'$ if necessary, we can assume that the homeomorphism is orientation-preserving. We deduce that there is an (orientation-preserving) homeomorphism from $Y'$ to some Dehn surgery $Y_{m/n}(K)$ that maps $K'$ with the dual knot. From now on we identify $Y'$ with $Y_{m/n}(K)$, and $K'$ to the dual knot.  By Lemma~\ref{lem:nullhomotopic-core}, we have that $Y$ and $Y'$  have isomorphic fundamental groups. The assumptions that $Y$ is irreducible and not a lens space imply (by geometrization) that $Y$ and $Y'$ are homeomorphic to each other; see for example \cite[Theorem 2.1.3]{3mBook}. 

We claim that the homeomorphism between $Y$ and $Y'$ can be chosen to be orientation-preserving. For this, we resort to Lemma~\ref{lem:nullhomotopic-core} again, which says that the surgery coefficient is $1/n$ for some integer $n$. If $n=0$, the surgery is trivial and $Y=Y'$. If $n \neq 0$,  as in the proof of Lemma~\ref{lem:nullhomotopic-core}, we have a cobordism $W$ from $Y$ to $Y'$ consisting of several $2$-handles, such that the map $\pi_1(Y)\to \pi_1(W)$ induced by inclusion is an isomorphism. In particular, the attaching circles are nullhomotopic. By extending these nullhomotopies to the 2-handles themselves, and then restricting to the outgoing boundary $Y'$, we obtain a degree-1 map $f: Y' \to Y$ that induces an isomorphism on $\pi_1$. 

We now distinguish two cases. If $\pi_1(Y)$ is infinite, then $Y$ is aspherical (see \cite[Section 3.2, (C.1)–
(C.2)]{3mBook}). By \cite[Theorem 2.1.2(2)]{3mBook}, there is a homeomorphism $h: Y' \to Y$ such that $h_*=f_*$ on $\pi_1$, so $\deg(h)=\deg(f)=1$ and $h$ is orientation-preserving. If $\pi_1(Y)$ is finite, $Y$ must be an elliptic 3-manifold. The hypothesis tells us it is not a lens space. This implies that $Y$ has no self-map of degree $-1$; see \cite[Example 2.6]{Sun}. We deduce that the homeomorphism we constructed from $Y$ to $Y'$ in the first place must have been orientation-preserving. 

In both cases, $Y$ and $Y'$ are orientation-preserving homeomorphic. We may assume that $Y=Y'$, and by applying Theorem~\ref{thm:Boileau-2} we conclude that $K =K'$.
\end{proof}

%!TEX root = main.tex

\section{An example illustrating why orientations matter}\label{sec:example}
As discussed in the introduction, there are pairs of knots in distinct lens spaces with orientation-preserving homeomorphic exteriors.  In this section, we show where the arguments used in this paper to prove Theorem~\ref{thm:Boileau} break down for such an example.  A key computation is Lemma~\ref{lem:index}, which shows that when $Y_{-1}(K)$ is orientation-preserving homeomorphic to $Y$, the index of a flat connection in $\chi^\sharp(W)$ associated to the 2-handle cobordism $W:Y \to Y_{-1}(K)$ is 0; that the two ends have the same orientation is used to be able to glue up the incoming and the outgoing boundary components of multiple copies of the cobordism $W$ to build a closed 4-manifold where we can use index formulas to compute the index of the ASD operator.  In this section, we carry out the index computation for the example from the introduction given in Figure~\ref{lens}. To get a $(-1)$-framed 2-handle cobordism, we reverse orientation and view $K$ as a knot in $-L(5,1)$ with a $-1$-surgery to $L(5,1)$.

\begin{figure}[ht!]
\labellist
\small\hair 2pt
\pinlabel {$\langle -5 \rangle$} at -35 246
\pinlabel {\textcolor{red}{$-1$}} [t] at 375 285
\endlabellist
\centering
\includegraphics[scale=.3]{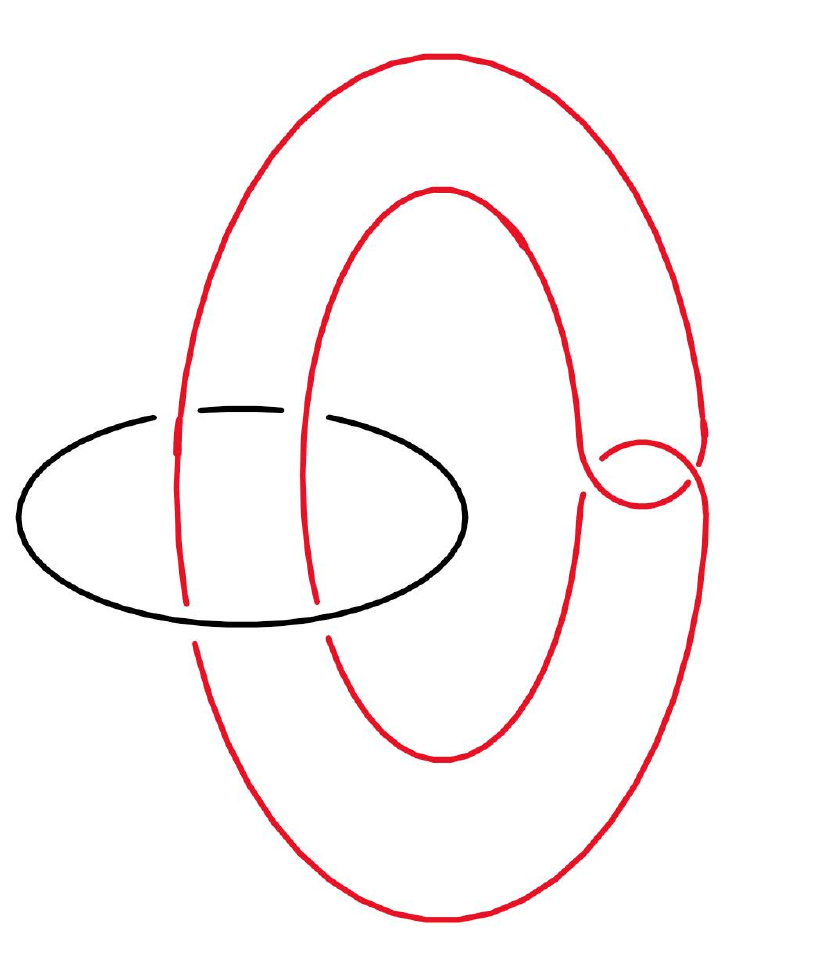}
\caption{The 2-handle cobordism $W$ from $-L(5,1)$ to $L(5,1)$.}\label{fig:W}
\end{figure}
\begin{prop}
	Let $K$ be the above knot in $-L(5,1)$ with $-1$-surgery to $L(5,1)$. Let \[W: -L(5,1) \to L(5,1)\] be the associated 2-handle cobordism and 
	$\rho:\pi_1(W)\to SU(2)$ be a non-trivial abelian representation. Then for the projectively flat connection $A^\sharp$ on $W^\sharp$ associated to $\rho$, the index of $\mathcal{D}_{A^\sharp}$ is non-zero.  
\end{prop}

In particular, this proposition shows that Lemma \ref{lem:index} fails if we relax the orientation-preserving homeomorphism assumption to an arbitrary homeomorphism.

\begin{proof}
Throughout this proof, we will write $Y$ for the lens space $L(5,1)$. In particular, $W^\sharp$ is a cobordism from $-Y^\sharp$ to $Y^\sharp$. In the same way as in the proof of Lemma \ref{lem:index}, we see that the kernel of the operator $\mathcal{D}_{A^\sharp}$ is trivial. To verify that the cokernel of $\mathcal{D}_{A^\sharp}$ is non-zero, we show that $H^+(W^\sharp; {\rm ad}_{A^\sharp})$ is non-trivial.  By excision, it suffices to show $H^+(W; {\rm ad}_{A})$ is non-zero, where $A$ is the flat abelian connection on $W$ associated to $\rho$. The local system ${\rm ad}_A$ splits as $\underline{\mathbb{R}}\oplus \underline{\Bbb C}_\rho$ where $\underline{\mathbb{R}}$ is the trivial bundle over $W$ and $\underline{\Bbb C}_\rho$ is a flat complex line bundle determined by $\rho$. More specifically, we may assume that $\rho$ maps the generator of $\pi_1(W)$ to the diagonal matrix in $SU(2)$ with diagonal entries $\omega^k$ and $\omega^{-k}$ where $\omega=e^{2\pi i/5}$ and $1\leq k\leq 4$ is an integer. Then $\underline{\Bbb C}_\rho$ is the quotient of the trivial line bundle $\underline{\Bbb C}$ over the universal cover $\widetilde W$ of $W$ by the action of $\pi_1(W)$ which maps the fiber over $x\in \widetilde W$ to the fiber over $\tau(x)$ by multiplication by $\omega^{2k}$. Here $\tau:\widetilde W \to \widetilde W$ is the covering action corresponding to the generator of $\pi_1(W)$.
Since $W$ is negative-definite, we have the splitting
\begin{align*}
	H^+(W; {\rm ad}_{A}) &\cong H^+(W;\underline{\mathbb{R}}) \oplus H^+(W;\underline{\Bbb C}_\rho) \\
	& \cong H^+(W;\underline{\Bbb C}_\rho)
\end{align*}
Therefore, we would like to see that $H^+(W;\underline{\Bbb C}_\rho)$ is non-zero. 

A relative Kirby diagram for $W$ is given in Figure~\ref{fig:W}, and from this we may see that Figure~\ref{fig:Wtilde} gives a relative Kirby diagram for the universal cover $\widetilde{W}: -\widetilde{Y} \to \widetilde{Y}$. Here $\widetilde{Y}$, the universal cover of $Y$, is the 3-dimensional sphere. 
\begin{figure}[ht!]
\labellist
\small\hair 2pt
\pinlabel {$\langle -1 \rangle$} at -25 170
\pinlabel {\textcolor{red}{$1$}} [t] at 435 330 
\pinlabel {\textcolor{red}{$1$}} [t] at 75 105
\pinlabel {\textcolor{red}{$1$}} [t] at 150 20
\pinlabel {\textcolor{red}{$1$}} [t] at 150 380
\pinlabel {\textcolor{red}{$1$}} [t] at 425 95
\endlabellist
\centering
\includegraphics[scale=.3]{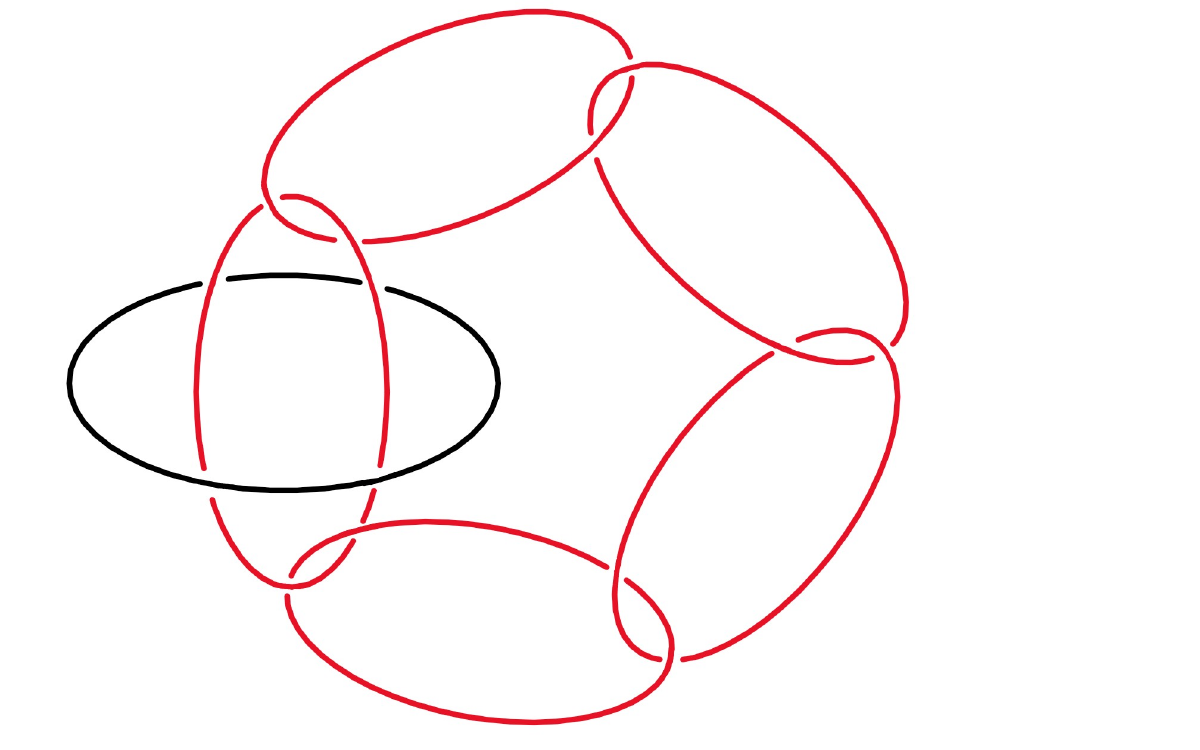}
\caption{The universal cover, $\widetilde{W}$, of the 2-handle cobordism $W$.}\label{fig:Wtilde}
\end{figure}
A basis $\{v_1,v_2,\dots,v_5\}$ for the second cohomology group of $\widetilde W$ is given by the Poincar\'e duals of the red 2-handles in Figure~\ref{fig:Wtilde}. The induced action by the generator of $\pi_1(W) =\Bbb Z/5$ is given by the obvious cyclic permutation $v_i \to v_{i+1}$. If we orient the red curves in Figure~\ref{fig:Wtilde} in a $\Bbb Z/5$-invariant fashion, each connected component of the red link has linking number $-1$ with the neighboring components and linking number $0$ with the remaining two components. By the $\Bbb Z/5$-symmetry the framings of the red components are equal to the same integer $a$. Using the computation of the linking numbers, we see that the sum $\widetilde \zeta$ of the cohomology classes $v_i$ satisfies the self-intersection relation
\[
  \widetilde \zeta \cdot \widetilde \zeta=5(a-2).
\] 
Since $\widetilde \zeta$ is the pullback of the cohomology class $\zeta$ associated to the red curve in Figure~\ref{fig:W}, we have
\[
  \widetilde \zeta \cdot \widetilde \zeta=5\zeta \cdot \zeta=-5.
\]
From this, we conclude that the framings of the red curves in Figure~\ref{fig:Wtilde} equal $a=1$, as it is indicated in the figure. In particular, the intersection form of $\widetilde{W}$, with respect to the above basis, is  
\begin{equation}\label{int-form}
Q_{\widetilde{W}} = \left[\begin{array}{ccccc} 1 & -1 & 0 & 0 & -1 \\ -1 & 1 & -1 & 0 & 0 \\ 0 & -1 & 1 & -1 & 0 \\ 0 & 0 & -1 & 1 & -1 \\ -1 & 0 & 0 & -1 & 1  \end{array}\right]. 
\end{equation}
In the above presentation, the $i^{\rm th}$ row and the $i^{\rm th}$ column correspond to the basis element $v_i$.

Now, we can see easily that $H^2(W;\underline{\Bbb C}_\rho)$ is the 1-dimensional complex vector space generated by $\sum_{i} \omega^{2ik}v_i$, and the Hermitian form on $H^2(W;\underline{\Bbb C}_\rho)$ determined by \eqref{int-form} is positive definite on this vector space. In summary, the cokernel of $\mathcal{D}_{A^\sharp}$ is $2$-dimensional. 
\end{proof}

\bibliography{references}
\bibliographystyle{alpha.bst}
\end{document}